\documentclass{gtart}

\def\ifplaintex{\expandafter\ifx\csname documentclass\endcsname\relax}

\ifplaintex 
\else
\fi

\expandafter\ifx\csname beginpicture\endcsname\relax
\expandafter\ifx\csname documentclass\endcsname\relax
\input pictex \else
\input prepictex \input pictex \input postpictex \fi\fi

\def\gt{{\mathsurround=0pt\it $\cal G\mskip-2mu$eometry \&\ 
$\cal T\!\!$opology}}        

\def\gtp{{\mathsurround=0pt\it $\cal G\mskip-2mu$eometry \&\ 
$\cal T\!\!$opology $\cal P\!$ublications}}  

\def\lognumber#1{\def\thelognumber{#1}}
\def\volumenumber#1{\def\thevolumenumber{#1}}
\def\papernumber#1{\def\thepapernumber{#1}}
\def\volumeyear#1{\def\thevolumeyear{#1}}

\def\pagenumbers#1#2{\def\startpage{#1}\def\finishpage{#2}}
\def\published#1{\def\publishdate{#1}}
\def\proposed#1{\def\theproposer{#1}}
\def\seconded#1{\def\theseconders{#1}}
\def\received#1{\def\receiveddate{#1}}
\def\revised#1{\def\reviseddate{#1}}
\def\accepted#1{\def\accepteddate{#1}}

\long\def\asciiabstract#1{\long\def\theasciiabstract{#1}}
\def\asciikeywords#1{\def\theasciikeywords{#1}}

\let\\\par\let\thelognumber\relax
\let\thevolumenumber\relax\let\thepapernumber\relax
\let\thevolumeyear\relax\let\thesamplenumber\relax\let\startpage\relax
\let\finishpage\relax\let\publishdate\relax\let\receiveddate\relax
\let\reviseddate\relax\let\accepteddate\relax\let\theasciititle\relax
\let\theasciiauthors\relax
\let\theasciiabstract\relax\let\theasciikeywords\relax
\let\theasciiemail\relax\let\theshortauthors\relax\let\theshorttitle\relax

\long\def\maketitlep{   

\count0=\startpage

\gt\hfill      
\beginpicture
\setcoordinatesystem units <0.33truein, 0.33truein> point at 2.2 0.9
\setplotsymbol ({$\cal G$})
\plotsymbolspacing=9truept
\circulararc 315 degrees from 0 1 center at 0 0
\setplotsymbol ({$\cal T$})
\circulararc 315 degrees from 1 -1 center at 1 0
\endpicture
\break
{\small\ifx\thesamplenumber\relax 
Volume \else Sample
\fi\thevolumenumber\ (\thevolumeyear)
\startpage--\finishpage\nl
Published: \publishdate}
\vglue 0.5truein plus 0.4fil minus 0.1truein

{\parskip=0pt\leftskip 0pt plus 1fil\def\\{\par\smallskip}{\ifplaintex\large
\else\Large\fi\bf\thetitle}\par\medskip}   

\vglue 0pt plus 0.1fil 

{\parskip=0pt\leftskip 0pt plus 1fil\def\\{\par}{\sc\theauthors}
\par\medskip}

\vglue 0pt plus 0.1fil 

{\small\parskip=0pt\let\newline\\
{\leftskip 0pt plus 1fil\def\\{\par}{\sl\theaddress}\par}
\expandafter\ifx\theemail\relax    
\relax\else\vglue 5pt plus 0.02fil minus 2pt\def\\{\stdspace{\rm 
and}\stdspace} 
\cl{Email:\stdspace\tt\theemail}\fi
\ifx\theurl\relax                  
\relax\else\vglue 5pt plus 0.02fil minus 2pt\def\\{\stdspace{\rm 
and}\stdspace}
\cl{URL:\stdspace\tt\theurl}\fi\par}

\vglue 7pt plus 0.3fil minus 3pt

{\bf Abstract}
\vglue 5pt plus 0.1fil minus 2pt

\theabstract

\vglue 7pt plus 0.3fil minus 3pt

{\bf AMS Classification numbers}\quad Primary:\quad \theprimaryclass

Secondary:\quad \thesecondaryclass

\vglue 5pt plus 0.3fil minus 2pt

{\bf Keywords:}\quad \thekeywords

\vglue 10pt plus 0.5fil minus 5pt

{\small  Proposed: \theproposer\hfill Received: \receiveddate\nl
Seconded: \theseconders\hfill 
\ifx\reviseddate\relax                         
Accepted: \accepteddate                        
\else
Revised: \reviseddate                          
\fi}
\eject
}       

\let\maketitlepage\maketitlep
\let\maketitle\maketitlepage

\font\phead=cmsl9 scaled 950
\font=cmsl9 scaled 1050
\font\pnum=cmbx10 scaled 913
\font\lnum=cmbx10 
\font\pfoot=cmsl9 scaled 950
\font=cmsl9 scaled 1050
\ifplaintex
\headline{\vbox to 0pt{\vskip -4.5mm\line{\small\phead\ifnum
\count0=\startpage ISSN 1364-0380 (on line)
1465-3060 (printed) \hfill {\pnum\folio}\else\ifodd\count0\def\\{ }%
\ifx\theshorttitle\relax\thetitle\else\theshorttitle\fi\hfill{\pnum\folio}
\else\def\\{ and }{\pnum\folio}\hfill\ifx\theshortauthors\relax\theauthors
\else\theshortauthors\fi\fi\fi}\vss}}
\footline{\vbox to 0pt{\vglue 0mm\line{\small\pfoot\ifnum\count0=\startpage
\copyright\ \gtp\hfill\else
\gt, Volume \thevolumenumber\ (\thevolumeyear)\hfill\fi}\vss
}}
\else
\makeatletter
\def\@oddhead{{\small\lhead\ifnum\count0=\startpage ISSN 1364-0380 (on line)
1465-3060 (printed) \hfill {\lnum\number\count0}\else\ifodd\count0
\def\\{ }\ifx\theshorttitle\relax \thetitle \else\theshorttitle\fi\hfill
{\lnum\number\count0}\else\def\\{ and }{\lnum\number\count0}
\hfill\ifx\theshortauthors\relax 
\theauthors\else\theshortauthors\fi\fi\fi}}\def\@evenhead{\@oddhead}
\def\@oddfoot{\small\lfoot\ifnum\count0=\startpage\copyright\ \gtp\hfill\else
\gt, Volume \thevolumenumber\ (\thevolumeyear)\hfill\fi}
\def\@evenfoot{\@oddfoot}
\makeatother
\fi

\newwrite\gtoutfile
\long\gdef\makeheadfile{  
{\def\\{, }\def\s{ }
\immediate\openout\gtoutfile head.xxx
\immediate\write\gtoutfile{To: math@arxiv.org}
\immediate\write\gtoutfile{Subject: put or rep NNNNN:pppp}
\immediate\write\gtoutfile{--text follows this line--}
\immediate\write\gtoutfile{Proxy-for: \ifx\theasciiauthors\relax
\theauthors\else\theasciiauthors\fi\s<\ifx\theasciiemail\relax\theemail\else\theasciiemail\fi>}
\immediate\write\gtoutfile{\noexpand\\}
\immediate\write\gtoutfile{Authors: \ifx\theasciiauthors\relax
\theauthors\else\theasciiauthors\fi}
\immediate\write\gtoutfile{Title: \ifx\theasciititle\relax
\thetitle\else\theasciititle\fi}
\immediate\write\gtoutfile{Subj-class: GT or SG or MG etc}
\immediate\write\gtoutfile{MSC-class: \theprimaryclass\ifx\thesecondaryclass\relax\else, \thesecondaryclass\fi}
\immediate\write\gtoutfile{Journal-ref: Geom. Topol. \thevolumenumber
(\thevolumeyear) \startpage-\finishpage}
\immediate\write\gtoutfile{Comments: Published by Geometry and Topology at}
\immediate\write\gtoutfile{\s\s http://www.maths.warwick.ac.uk/gt/GTVol\thevolumenumber/paper\thepapernumber.abs.html}
\immediate\write\gtoutfile{\noexpand\\}
\immediate\write\gtoutfile{}
\ifx\theasciiabstract\relax
\immediate\write\gtoutfile{\theabstract}\else
\immediate\write\gtoutfile{\theasciiabstract}\fi
\immediate\write\gtoutfile{}
\immediate\write\gtoutfile{\noexpand\\}
\immediate\write\gtoutfile{}
\immediate\closeout\gtoutfile}}  

\def\maketitlepage{\maketitlep\makeheadfile}
\let\maketitle\maketitlepage

\def\ifplaintex{\expandafter\ifx\csname documentclass\endcsname\relax}

\ifplaintex 
\else
\fi

\expandafter\ifx\csname beginpicture\endcsname\relax
\expandafter\ifx\csname documentclass\endcsname\relax
\input pictex \else
\input prepictex \input pictex \input postpictex \fi\fi

\def\gt{{\mathsurround=0pt\it $\cal G\mskip-2mu$eometry \&\ 
$\cal T\!\!$opology}}        

\def\gtp{{\mathsurround=0pt\it $\cal G\mskip-2mu$eometry \&\ 
$\cal T\!\!$opology $\cal P\!$ublications}}  

\def\lognumber#1{\def\thelognumber{#1}}
\def\volumenumber#1{\def\thevolumenumber{#1}}
\def\papernumber#1{\def\thepapernumber{#1}}
\def\volumeyear#1{\def\thevolumeyear{#1}}

\def\pagenumbers#1#2{\def\startpage{#1}\def\finishpage{#2}}
\def\published#1{\def\publishdate{#1}}
\def\proposed#1{\def\theproposer{#1}}
\def\seconded#1{\def\theseconders{#1}}
\def\received#1{\def\receiveddate{#1}}
\def\revised#1{\def\reviseddate{#1}}
\def\accepted#1{\def\accepteddate{#1}}

\long\def\asciiabstract#1{\long\def\theasciiabstract{#1}}
\def\asciikeywords#1{\def\theasciikeywords{#1}}

\let\\\par\let\thelognumber\relax
\let\thevolumenumber\relax\let\thepapernumber\relax
\let\thevolumeyear\relax\let\thesamplenumber\relax\let\startpage\relax
\let\finishpage\relax\let\publishdate\relax\let\receiveddate\relax
\let\reviseddate\relax\let\accepteddate\relax\let\theasciititle\relax
\let\theasciiauthors\relax
\let\theasciiabstract\relax\let\theasciikeywords\relax
\let\theasciiemail\relax\let\theshortauthors\relax\let\theshorttitle\relax

\long\def\maketitlep{   

\count0=\startpage

\gt\hfill      
\beginpicture
\setcoordinatesystem units <0.33truein, 0.33truein> point at 2.2 0.9
\setplotsymbol ({$\cal G$})
\plotsymbolspacing=9truept
\circulararc 315 degrees from 0 1 center at 0 0
\setplotsymbol ({$\cal T$})
\circulararc 315 degrees from 1 -1 center at 1 0
\endpicture
\break
{\small\ifx\thesamplenumber\relax 
Volume \else Sample
\fi\thevolumenumber\ (\thevolumeyear)
\startpage--\finishpage\nl
Published: \publishdate}
\vglue 0.5truein plus 0.4fil minus 0.1truein

{\parskip=0pt\leftskip 0pt plus 1fil\def\\{\par\smallskip}{\ifplaintex\large
\else\Large\fi\bf\thetitle}\par\medskip}   

\vglue 0pt plus 0.1fil 

{\parskip=0pt\leftskip 0pt plus 1fil\def\\{\par}{\sc\theauthors}
\par\medskip}

\vglue 0pt plus 0.1fil 

{\small\parskip=0pt\let\newline\\
{\leftskip 0pt plus 1fil\def\\{\par}{\sl\theaddress}\par}
\expandafter\ifx\theemail\relax    
\relax\else\vglue 5pt plus 0.02fil minus 2pt\def\\{\stdspace{\rm 
and}\stdspace} 
\cl{Email:\stdspace\tt\theemail}\fi
\ifx\theurl\relax                  
\relax\else\vglue 5pt plus 0.02fil minus 2pt\def\\{\stdspace{\rm 
and}\stdspace}
\cl{URL:\stdspace\tt\theurl}\fi\par}

\vglue 7pt plus 0.3fil minus 3pt

{\bf Abstract}
\vglue 5pt plus 0.1fil minus 2pt

\theabstract

\vglue 7pt plus 0.3fil minus 3pt

{\bf AMS Classification numbers}\quad Primary:\quad \theprimaryclass

Secondary:\quad \thesecondaryclass

\vglue 5pt plus 0.3fil minus 2pt

{\bf Keywords:}\quad \thekeywords

\vglue 10pt plus 0.5fil minus 5pt

{\small  Proposed: \theproposer\hfill Received: \receiveddate\nl
Seconded: \theseconders\hfill 
\ifx\reviseddate\relax                         
Accepted: \accepteddate                        
\else
Revised: \reviseddate                          
\fi}
\eject
}       

\let\maketitlepage\maketitlep
\let\maketitle\maketitlepage

\font\phead=cmsl9 scaled 950
\font=cmsl9 scaled 1050
\font\pnum=cmbx10 scaled 913
\font\lnum=cmbx10 
\font\pfoot=cmsl9 scaled 950
\font=cmsl9 scaled 1050
\ifplaintex
\headline{\vbox to 0pt{\vskip -4.5mm\line{\small\phead\ifnum
\count0=\startpage ISSN 1364-0380 (on line)
1465-3060 (printed) \hfill {\pnum\folio}\else\ifodd\count0\def\\{ }%
\ifx\theshorttitle\relax\thetitle\else\theshorttitle\fi\hfill{\pnum\folio}
\else\def\\{ and }{\pnum\folio}\hfill\ifx\theshortauthors\relax\theauthors
\else\theshortauthors\fi\fi\fi}\vss}}
\footline{\vbox to 0pt{\vglue 0mm\line{\small\pfoot\ifnum\count0=\startpage
\copyright\ \gtp\hfill\else
\gt, Volume \thevolumenumber\ (\thevolumeyear)\hfill\fi}\vss
}}
\else
\makeatletter
\def\@oddhead{{\small\lhead\ifnum\count0=\startpage ISSN 1364-0380 (on line)
1465-3060 (printed) \hfill {\lnum\number\count0}\else\ifodd\count0
\def\\{ }\ifx\theshorttitle\relax \thetitle \else\theshorttitle\fi\hfill
{\lnum\number\count0}\else\def\\{ and }{\lnum\number\count0}
\hfill\ifx\theshortauthors\relax 
\theauthors\else\theshortauthors\fi\fi\fi}}\def\@evenhead{\@oddhead}
\def\@oddfoot{\small\lfoot\ifnum\count0=\startpage\copyright\ \gtp\hfill\else
\gt, Volume \thevolumenumber\ (\thevolumeyear)\hfill\fi}
\def\@evenfoot{\@oddfoot}
\makeatother
\fi

\newwrite\gtoutfile
\long\gdef\makeheadfile{  
{\def\\{, }\def\s{ }
\immediate\openout\gtoutfile head.xxx
\immediate\write\gtoutfile{To: math@arxiv.org}
\immediate\write\gtoutfile{Subject: put or rep NNNNN:pppp}
\immediate\write\gtoutfile{--text follows this line--}
\immediate\write\gtoutfile{Proxy-for: \ifx\theasciiauthors\relax
\theauthors\else\theasciiauthors\fi\s<\ifx\theasciiemail\relax\theemail\else\theasciiemail\fi>}
\immediate\write\gtoutfile{\noexpand\\}
\immediate\write\gtoutfile{Authors: \ifx\theasciiauthors\relax
\theauthors\else\theasciiauthors\fi}
\immediate\write\gtoutfile{Title: \ifx\theasciititle\relax
\thetitle\else\theasciititle\fi}
\immediate\write\gtoutfile{Subj-class: GT or SG or MG etc}
\immediate\write\gtoutfile{MSC-class: \theprimaryclass\ifx\thesecondaryclass\relax\else, \thesecondaryclass\fi}
\immediate\write\gtoutfile{Journal-ref: Geom. Topol. \thevolumenumber
(\thevolumeyear) \startpage-\finishpage}
\immediate\write\gtoutfile{Comments: Published by Geometry and Topology at}
\immediate\write\gtoutfile{\s\s http://www.maths.warwick.ac.uk/gt/GTVol\thevolumenumber/paper\thepapernumber.abs.html}
\immediate\write\gtoutfile{\noexpand\\}
\immediate\write\gtoutfile{}
\ifx\theasciiabstract\relax
\immediate\write\gtoutfile{\theabstract}\else
\immediate\write\gtoutfile{\theasciiabstract}\fi
\immediate\write\gtoutfile{}
\immediate\write\gtoutfile{\noexpand\\}
\immediate\write\gtoutfile{}
\immediate\closeout\gtoutfile}}  

\def\maketitlepage{\maketitlep\makeheadfile}
\let\maketitle\maketitlepage

\lognumber{88}
\volumenumber{4}\papernumber{19}\volumeyear{2000}
\pagenumbers{537}{579}
\proposed{Steve Ferry}
\seconded{Robion Kirby, Cameron Gordon}
\received{30 July 1999}
\revised{8 December 2000}
\accepted{27 December 2000}
\published{27 December 2000}

\usepackage{graphicx,amsmath,amssymb}

\theoremstyle{plain}
\newtheorem{theorem}{Theorem}
\newtheorem{claim}{Claim}
\newtheorem{corollary}{Corollary}
\newtheorem{lemma}{Lemma}
\newtheorem{proposition}{Proposition}

\theoremstyle{remark}
\newtheorem{remark}{Remark}
\newtheorem{case}{Case}
\newtheorem{example}{Example}

\numberwithin{equation}{section}
\nocolon  

\begin{document}

\title{Manifolds with non-stable fundamental groups\\at infinity}
\author{Craig R Guilbault}
\address{Department of Mathematical Sciences\\
University of Wisconsin-Milwaukee\\WI 53201-0413, USA}
\email{craigg@uwm.edu}

\primaryclass{57N15, 57Q12}
\secondaryclass{57R65, 57Q10}

\keywords{Non-compact manifold, ends, collar, homotopy collar, pseudo-collar,
semistable, Mittag--Leffler}

\asciikeywords{Non-compact manifold, ends, collar, homotopy collar,
pseudo-collar, semistable, Mittag-Leffler}

\begin{abstract}
The notion of an open collar is generalized to that of a pseudo-collar.
Important properties and examples are discussed. The main result gives
conditions which guarantee the existence of a pseudo-collar structure on the
end of an open $n$--manifold ($n\geq7$). This paper may be viewed as a
generalization of Siebenmann's famous collaring theorem to open manifolds with
non-stable fundamental group systems at infinity.
\end{abstract}
\asciiabstract{
The notion of an open collar is generalized to that of a pseudo-collar.
Important properties and examples are discussed. The main result gives
conditions which guarantee the existence of a pseudo-collar structure on the
end of an open n-manifold (n > 6). This paper may be viewed as a
generalization of Siebenmann's famous collaring theorem to open manifolds with
non-stable fundamental group systems at infinity.}

\maketitlepage

\section{Introduction}

One of the best known and most frequently applied theorems in the study of
non-compact manifolds is found in L\,C Siebenmann's 1965 PhD thesis. It
gives necessary and sufficient conditions for the end of an open manifold to
possess the simplest possible structure---that of an open collar.

\begin{theorem}[From \cite{Si1}]
\label{sieb}A one ended open $n$--manifold $M^{n}$ ($n\geq6
$) contains an open collar neighborhood of infinity if and only if each of the
following is satisfied:

\begin{enumerate}
\item $M^{n}$ is inward tame at infinity,

\item $\pi_{1}\left(  \varepsilon(M^{n})\right)  $ is stable, and

\item $\sigma_{\infty}\left(  M^{n}\right)  \in\widetilde{K}_{0}\left(
\mathbb{Z[}\pi_{1}(\varepsilon(M^{n}))]\right)  $ is trivial.
\end{enumerate}
\end{theorem}

A \emph{neighborhood of infinity }is a subset $U\subset M^{n}$ with the
property that $\overline{M^{n}-U}$ is compact. We say that $U$ is an
\emph{open collar} if it is a manifold with compact boundary and
$U\approx\partial U\times\lbrack0,\infty)$. Other terminology and notation
used in this theorem will be discussed later.

\begin{remark}
A $3$--dimensional version of Theorem \ref{sieb} may be found in \cite{HP},
while a $5$--dimensional version (with some restrictions) may be found in
\cite{FQ}. In \cite{KS}, it is shown that Theorem \ref{sieb} fails in
dimension $4$.
\end{remark}

One of the beauties of Theorem \ref{sieb} is the simple structure it places on
the ends of certain manifolds. At the same time, this simplicity greatly
limits the class of manifolds to which the theorem applies. Indeed, many
interesting and important non-compact manifolds are ``too complicated at
infinity'' to be collarable. Frequently the condition these manifolds violate
is $\pi_{1}$--stability. In this paper we present a program to generalize
Theorem \ref{sieb} so that it applies to manifolds with non-stable fundamental
groups at infinity. Of course, a manifold with non-stable fundamental group at
infinity cannot be collarable, so we must be satisfied with a less rigid
structure on its end. The structure we have chosen to pursue will be called a
\emph{pseudo-collar}.

We say that a manifold $U^{n}$ with compact boundary is a \emph{homotopy
collar} provided the inclusion $\partial U^{n}\hookrightarrow U^{n}$ is a
homotopy equivalence. As it turns out, a homotopy collar may possess very
little additional structure, hence, we define the following more rigid notion.
A \emph{pseudo-collar} is a homotopy collar that contains arbitrarily small
homotopy collar neighborhoods of infinity.

With the above definition established, the goal of this paper can be described
as a study of pseudo-collar\-abil\-ity in high dimensional manifolds. For the
sake of the experts, we state our principal result now. A more thorough
development and motivation of this theorem can be found in Section
\ref{pseudo}.\medskip

\noindent\textbf{Main Existence Theorem}\qua\emph{A one ended open }%
$n$\emph{--manifold }$M^{n}$\emph{ (}$n\geq7$\emph{) is pseudo-collarable
provided each of the following is satisfied:}

\begin{enumerate}
\item $M^{n}$\emph{ is inward tame at infinity,}

\item $\pi_{1}(\varepsilon(M^{n}))$\emph{ is perfectly semistable,}

\item $\sigma_{\infty}\left(  M^{n}\right)  =0\in\widetilde{K}_{0}\left(
\pi_{1}\left(  \varepsilon(M^{n}\right)  \right)  )$\emph{, and}

\item $\pi_{2}(\varepsilon(M^{n}))$\emph{ is semistable.\smallskip}
\end{enumerate}

As the reader can see, Condition 1 of Theorem \ref{sieb} is unchanged in our
more general setting. Condition 3 has been reformulated so that it applies to
situations where the fundamental group at infinity is not stable---but it also
is essentially unchanged. Both of these conditions are discussed in Section
\ref{background}. The weakening of Condition 2 is the main task in this
paper.\ Much of our work is done with no restrictions on the fundamental group
at infinity; however, it eventually becomes necessary to focus on manifolds
with \emph{perfectly semistable} $\pi_{1}$--systems at infinity. These systems
are semistable (also called \emph{Mittag--Leffler}) and have bonding maps with
perfect kernels. Condition 4 is different from the others---it has no analog
in Theorem \ref{sieb}, and we are not sure whether it is necessary. It does,
however, play a crucial role in our proof.

Semistability conditions are well-established in studies of non-compact
$3$--man\-ifolds (see \cite{Ka} or \cite{BT}) and also in studies of ends of
groups (see \cite{Mi}), so it seems fitting that they play a role in the study
of high-dimensional manifolds. Precise definitions of these conditions may be
found in Section \ref{groups}.

In Section \ref{background} we review some basics in the study of non-compact
manifolds, then in Section \ref{pseudo} we explore the topology of
pseudo-collars. Some examples are discussed and basic geometric and algebraic
properties are derived. These provide the necessary framework and motivation
for our Main Existence Theorem. Most of the remainder of the paper (Sections
\ref{outline}--\ref{n-2}) is geared towards proving this theorem. The strategy
is much the same as that used in \cite{Si1}; however, since the hypothesis of
$\pi_{1}$--stability is thoroughly ingrained in Siebenmann's work, nearly all
steps require some revision. Sometimes these revisions are significant, while
at other times the original arguments already suffice. For completeness,
portions of \cite{Si1} have been repeated. The reader who makes it to the end
of this paper will reprove Siebenmann's theorem in the process. (See Remark
\ref{siebproof}.) In the final section of this paper we discuss some open questions.

We conclude this introduction by defending our choice of
``pseudo-collar\-abil\-ity'' as the appropriate generalization of collarability.

At first glance, one might expect ``homotopy collar'' to be a good enough
generalization of ``collar''. Unfortunately, homotopy collars carry very
little useful structure beyond what is given by their definition. For example,
every contractible open manifold (no matter how badly behaved at infinity)
contains a homotopy collar neighborhood of infinity---just consider the
complement of a small open ball in the manifold. Hence, some additional
structure is desired. Propositions \ref{decomposition} and \ref{laminated} and
Theorem \ref{necessary} show that pseudo-collars do indeed carry a great deal
of additional structure.

A second reason for defining pseudo-collars as we have is to mimic a key
property possessed by genuine collars. In particular, a collar structure on
the end of a manifold guarantees the existence of arbitrarily small collar
neighborhoods of that end. Although this observation is trivial, it is
extremely important in applications. It seems that any useful generalization
of \ ``collar'' should have an analogous property.

A third factor which focused our attention on pseudo-collar\-abil\-ity was
work by Chapman and Siebenmann on $\mathcal{Z}$--compactifications of Hilbert
cube manifolds. Although they advertise their main result as an infinite
dimensional version of Theorem \ref{sieb}, it is really much more general. In
particular, it applies to Hilbert cube manifolds with non-collarable ends.
Their program can be broken into two parts. First they determine necessary and
sufficient conditions for a one ended Hilbert cube manifold $X$ to contain
arbitrarily small neighborhoods $U$ of infinity for which
$Bdry(U)\hookrightarrow U$ is a homotopy equivalence. (In our language, they
determine when $X$ is pseudo-collarable.) Next they combine the structure
supplied by the pseudo-collar with some powerful results from Hilbert cube
manifold theory to determine whether a $\mathcal{Z}$--compactification is
possible. It is natural to ask if their program can be carried out in finite
dimensions. In this paper we focus on the first part of that program. We
intend to address the issue of $\mathcal{Z}$--compactifiability for finite
dimensional manifolds and its relationship to pseudo-collar\-abil\-ity in a
later paper.

A final reason for the choices we have made lies with some key examples and
current research trends in topology. For instance, the exotic universal
covering spaces produced by M Davis in \cite{Da} are all pseudo-collarable
but not collarable. Variations on those examples were produced in \cite{DJ}
with the aid of $CAT(0)$ geometry---they are also pseudo-collarable. Moreover,
many of the basic conditions necessary for pseudo-collar\-a\-bil\-ity are
satisfied by all $CAT\left(  0\right)  $ manifolds and also by universal
covers of all aspherical manifolds with word hyperbolic or $CAT\left(
0\right)  $ fundamental groups. Thus, the collection of examples to which our
techniques might be applied appears quite rich.

We wish to thank Steve Ferry for directing us to \cite{Fe1} and for sharing a
copy of \cite{Fe2} which contains a clear and concise exposition of
Siebenmann's thesis.

\section{Inverse sequences and group theory\label{groups}}

Throughout this section all arrows denote homomorphisms, while arrows of the
type $\twoheadrightarrow$ or $\twoheadleftarrow$ denote surjections. The
symbol $\cong$ denotes isomorphisms.

Let
\[
G_{0}\overset{\lambda_{1}}{\longleftarrow}G_{1}\overset{\lambda_{2}%
}{\longleftarrow}G_{2}\overset{\lambda_{3}}{\longleftarrow}\cdots
\]
be an inverse sequence of groups and homomorphisms. A \emph{subsequence} of
$\left\{  G_{i},\lambda_{i}\right\}  $ is an inverse sequence of the form
\[
G_{i_{0}}\overset{\lambda_{i_{0}+1}\circ\cdots\circ\lambda_{i_{1}}%
}{\longleftarrow}G_{i_{1}}\overset{\lambda_{i_{1}+1}\circ\cdots\circ
\lambda_{i_{2}}}{\longleftarrow}G_{i_{2}}\overset{\lambda_{i_{2}+1}\circ
\cdots\circ\lambda_{i_{3}}}{\longleftarrow}\cdots.
\]
In the future we will denote a composition $\lambda_{i}\circ\cdots\circ
\lambda_{j}$ ($i\leq j$) by $\lambda_{i,j}$.

We say that sequences $\left\{  G_{i},\lambda_{i}\right\}  $ and $\left\{
H_{i},\mu_{i}\right\}  $ are \emph{pro-equivalent} if, after passing to
subsequences, there exists a commuting diagram:
\[%
\begin{array}
[c]{ccccccc}%
G_{i_{0}} & {\buildrel{\lambda_{i_{0}+1,i_{1}}}\over{\longleftarrow}} & G_{i_{1}} &
\overset{\lambda_{i_{1}+1,i_{1}}}{\longleftarrow} & G_{i_{2}} & \overset
{\lambda_{i_{2}+1,i_{3}}}{\longleftarrow} & \cdots\\
& \nwarrow\quad\swarrow &  & \nwarrow\quad\swarrow &  & \nwarrow\quad\swarrow
& \\
&  H_{j_{0}} & \overset{\mu_{i_{0}+1,i_{1}}}{\longleftarrow} & H_{j_{1}} &
\overset{\mu_{i_{1}+1,i_{1}}}{\longleftarrow} & H_{j_{2}} & \cdots
\end{array}
\]
Clearly an inverse sequence is pro-equivalent to any of its subsequences. To
avoid tedious notation, we often do not distinguish $\left\{  G_{i}%
,\lambda_{i}\right\}  $ from its subsequences. Instead we simply assume that
$\left\{  G_{i},\lambda_{i}\right\}  $ has the desired properties of a
preferred subsequence---often prefaced by the words ``after passing to a
subsequence and relabelling''.

The \emph{inverse limit }of a sequence $\left\{  G_{i},\lambda_{i}\right\}  $
is a subgroup of $\prod G_{i}$ defined by
\[
{\displaystyle\lim_{\textstyle\longleftarrow}}\left\{  G_{i},\lambda_{i}\right\}  =\left\{  \left.
\left(  g_{0},g_{1},g_{2},\cdots\right)  \in\prod_{i=0}^{\infty}G_{i}\right|
\lambda_{i}\left(  g_{i}\right)  =g_{i-1}\right\}  .
\]
Notice that for each $i$, there is a \emph{projection homomorphism}
$p_{i}\co {\displaystyle\lim_{\textstyle\longleftarrow}}\left\{  G_{i},\lambda_{i}\right\}  \rightarrow
G_{i}$. It is a standard fact that pro-equivalent inverse sequences have
isomorphic inverse limits.

An inverse sequence $\left\{  G_{i},\lambda_{i}\right\}  $ is \emph{stable} if
it is pro-equivalent to a constant sequence $\left\{  H,id\right\}  $. It is
easy to see that $\left\{  G_{i},\lambda_{i}\right\}  $ is stable if and only
if, after passing to a subsequence and relabelling, there is a commutative
diagram of the form:
\[%
\begin{array}
[c]{ccccccccc}%
G_{0} & \overset{\lambda_{1}}{\longleftarrow} & G_{1} & \overset{\lambda_{2}%
}{\longleftarrow} & G_{2} & \overset{\lambda_{3}}{\longleftarrow} & G_{3} &
\overset{\lambda_{4}}{\longleftarrow} & \cdots\\
& \nwarrow\quad\swarrow &  & \nwarrow\quad\swarrow &  & \nwarrow\quad\swarrow
&  &  & \\
&  im(\lambda_{1}) & \overset{\cong}{\longleftarrow} & im(\lambda_{2}) &
\overset{\cong}{\longleftarrow} & im(\lambda_{3}) & \overset{\cong
}{\longleftarrow} & \cdots &
\end{array}
\]
In this case $H\cong{\displaystyle\lim_{\textstyle\longleftarrow}}\left\{  G_{i},\lambda_{i}\right\}
\cong im(\lambda_{i})$ and each projection homomorphism takes ${\displaystyle\lim_{\textstyle\longleftarrow}}\left\{  G_{i},\lambda_{i}\right\}  $ isomorphically onto the
corresponding $im(\lambda_{i})$.

The sequence $\left\{  G_{i},\lambda_{i}\right\}  $ is \emph{semistable} (or
\emph{Mittag--Leffler}) if it is pro-equivalent to an inverse sequence
$\left\{  H_{i},\mu_{i}\right\}  $ for which each $\mu_{i}$ is surjective.
Equivalently, $\left\{  G_{i},\lambda_{i}\right\}  $ is semistable if, after
passing to a subsequence and relabelling, there is a commutative diagram of
the form:
\[%
\begin{array}
[c]{ccccccccc}%
G_{0} & \overset{\lambda_{1}}{\longleftarrow} & G_{1} & \overset{\lambda_{2}%
}{\longleftarrow} & G_{2} & \overset{\lambda_{3}}{\longleftarrow} & G_{3} &
\overset{\lambda_{4}}{\longleftarrow} & \cdots\\
& \nwarrow\quad\swarrow &  & \nwarrow\quad\swarrow &  & \nwarrow\quad\swarrow
&  &  & \\
&  im(\lambda_{1}) & \twoheadleftarrow &  im(\lambda_{2}) & \twoheadleftarrow
&  im(\lambda_{3}) & \twoheadleftarrow & \cdots &
\end{array}
\]
We now describe a subclass of semistable inverse sequences which are of
particular interest to us. Recall that a \emph{commutator} element of a group
$H$ is an element of the form $xyx^{-1}y^{-1}$ where $x,y\in H$; and the
\emph{commutator subgroup} of $H,$ denoted $\left[  H,H\right]  $, is the
subgroup generated by all of its commutators. We say that $H$ is
\emph{perfect} if $\left[  H,H\right]  =H$. An inverse sequence of groups is
\emph{perfectly semistable} if it is pro-equivalent to an inverse sequence
\[
G_{0}\overset{\lambda_{1}}{\twoheadleftarrow}G_{1}\overset{\lambda_{2}%
}{\twoheadleftarrow}G_{2}\overset{\lambda_{3}}{\twoheadleftarrow}\cdots
\]
of finitely presentable groups and surjections where each $\ker\left(
\lambda_{i}\right)  $ \ is perfect. The following shows that inverse sequences
of this type behave well under passage to subsequences.

\begin{lemma}
\label{gp1}Suppose $f\co A\rightarrow B$ and $g\co B\rightarrow C$ are each
surjective group homomorphisms with perfect kernels. Then $g\circ
f\co A\rightarrow C$ is surjective and has perfect kernel.
\end{lemma}

\begin{proof}
Surjectivity is obvious. To see that $\ker\left(  g\circ f\right)  $ is
perfect, begin with $a\in A$ such that $\left(  g\circ f\right)  \left(
a\right)  =1$. Then $f(a)\in\ker(g)$, so by hypothesis we may write
\[
f\left(  a\right)  =\prod_{i=1}^{k}x_{i}y_{i}x_{i}^{-1}y_{i}^{-1}\text{ where
}x_{i},y_{i}\in\ker(g)\text{ for }i=1,\cdots,k.
\]
For each $i$, choose $u_{i},v_{i}\in A$ such that $f\left(  u_{i}\right)
=x_{i}$ and $f\left(  v_{i}\right)  =y_{i}$. Note that each $u_{i}$ and
$v_{i}$ lies in $\ker\left(  g\circ f\right)  $, and let
\[
a^{\prime}=\prod_{i=1}^{k}u_{i}v_{i}u_{i}^{-1}v_{i}^{-1}.
\]
Then $f\left(  a^{\prime}\right)  =f\left(  a\right)  $, which implies that
$a\left(  a^{\prime}\right)  ^{-1}\in\ker\left(  f\right)  $; so by hypothesis
we may write
\[
a\left(  a^{\prime}\right)  ^{-1}=\prod_{j=1}^{k}r_{j}s_{j}r_{j}^{-1}%
s_{j}^{-1}\text{ where }r_{j},s_{j}\in\ker(f)\text{ for }j=1,\cdots,l.
\]
Moreover, since $\ker(f)\subset\ker(g\circ f)$, each $r_{j},s_{j}$ lies in
$\ker(g\circ f).$

Finally, we write
\[
a=\left(  a\left(  a^{\prime}\right)  ^{-1}\right)  \cdot a^{\prime}=\left(
\prod_{j=1}^{k}r_{j}s_{j}r_{j}^{-1}s_{j}^{-1}\right)  \cdot\left(  \prod
_{i=1}^{k}u_{i}v_{i}u_{i}^{-1}v_{i}^{-1}\right)  ,
\]
which shows that $a\in\left[  \ker(g\circ f),\ker(g\circ f)\right]  $.
\end{proof}

\begin{corollary}
\label{gp2}If $\left\{  G_{i},\lambda_{i}\right\}  $ is an inverse sequence of
groups and surjections with perfect kernels, then so is any subsequence.
\end{corollary}

We conclude this section with three more group theoretic lemmas which will be
used later. The first is from \cite{Wa1}.

\begin{lemma}
\label{gp3}Let $A$ be a finitely generated group and $f\co A\rightarrow B$ and
$g\co B\rightarrow A$ be group homomorphisms with $f\circ g=id_{B}$. Then
$\ker(f)$ is the normal closure of a finite set of elements. Therefore, if $A$
is finitely presentable, then so is $B$.
\end{lemma}

\begin{proof}
Let $\left\{  a_{i}\right\}  _{i=1}^{k}$ be a generating set for $A$ and let
$X=\left\{  a_{i}\cdot(g\circ f)(a_{i}^{-1})\right\}  _{i=1}^{k}$. We will
show that $\ker(f)$ is the normal closure of $X$.

First note that $f\left(  a_{i}\cdot(g\circ f)(a_{i}^{-1})\right)  =f\left(
a_{i}\right)  \cdot\left(  f\circ g\circ f\right)  \left(  a_{i}^{-1}\right)
=f\left(  a_{i}\right)  \cdot f\left(  a_{i}^{-1}\right)  =1$, for each $i$.
Hence $X$ (and therefore the normal closure of $X)$, is contained in
$\ker\left(  f\right)  $.

As preparation for obtaining the reverse inclusion, let $w=w_{1}w_{2}\in
\ker\left(  f\right)  $ and observe that
\[%
\begin{array}
[c]{ll}%
w_{1}w_{2} & =w_{1}w_{2}\cdot(g\circ f)((w_{1}w_{2})^{-1})\\
& =w_{1}w_{2}\cdot(g\circ f)(w_{2}^{-1})\cdot(g\circ f)(w_{1}^{-1})\\
& =(w_{1}[w_{2}\cdot(g\circ f)(w_{2}^{-1})]w_{1}^{-1})(w_{1}\cdot(g\circ
f)(w_{1}^{-1})).
\end{array}
\]
With this identity as the main tool, induction on word length shows that
$\ker\left(  f\right)  $ $\subset normal$ $closure\left(  X\right)  $.
\end{proof}

The next lemma is from \cite{Si1}, where it is used for purposes similar to
our own.

\begin{lemma}
\label{gp4}Let $f\co A\twoheadrightarrow B$ be a group homomorphism, and suppose
$A=\left\langle a\mid r\right\rangle $ and $B=\left\langle b\mid
s\right\rangle $ are presentations with $\left|  a\right|  $ generators and
$\left|  s\right|  $ relators, respectively. Then $\ker(f)$ is the normal
closure of a set containing $\left|  a\right|  +\left|  s\right|  $ elements.
\end{lemma}

\begin{proof}
Let $\xi$ be a set of words so that $f\left(  a\right)  =\xi\left(  b\right)
$ in $B$. Since $f$ is surjective, there exists a set of words $\eta$ so that
$b=\eta\left(  f\left(  a\right)  \right)  $ in $B$. Then Tietze
transformations give the following isomorphisms:%
\begin{align*}
\left\langle b\mid s\right\rangle  &  \cong\left\langle a,b\mid a=\xi\left(
b\right)  ,s\left(  b\right)  \right\rangle \\
&  \cong\left\langle a,b\mid a=\xi\left(  b\right)  ,s\left(  b\right)
,r\left(  a\right)  ,b=\eta\left(  a\right)  \right\rangle \\
&  \cong\left\langle a,b\mid a=\xi\left(  \eta\left(  a\right)  \right)
,s\left(  \eta\left(  a\right)  \right)  ,r\left(  a\right)  ,b=\eta\left(
a\right)  \right\rangle \\
&  \cong\left\langle a\mid a=\xi\left(  \eta\left(  a\right)  \right)
,s\left(  \eta\left(  a\right)  \right)  ,r\left(  a\right)  \right\rangle
\end{align*}
Now $f$ is specified by the last presentation via the correspondence
$a\longmapsto a$. Hence $\ker(f)$ is the normal closure of the $\left|
a\right|  +\left|  s\right|  $ elements of $\xi\left(  \eta\left(  a\right)
\right)  $ and $s\left(  \eta\left(  a\right)  \right)  $.
\end{proof}

The following lemma was extracted from the proof of Theorem 4 in \cite{Fe2}.

\begin{lemma}
\label{gp5}Each semistable inverse sequence $\left\{  G_{i},\lambda
_{i}\right\} $ of finitely presented\break groups is pro-equivalent to
an inverse sequence $\left\{ G_{i}^{\prime},\mu_{i}\right\} $ of
finitely presented groups with surjective bonding maps.
\end{lemma}

\begin{proof}
After passing to a subsequence and relabelling we have a diagram:
\[%
\begin{array}
[c]{ccccccccc}%
G_{0} & \overset{\lambda_{1}}{\longleftarrow} & G_{1} & \overset{\lambda_{2}%
}{\longleftarrow} & G_{2} & \overset{\lambda_{3}}{\longleftarrow} & G_{3} &
\overset{\lambda_{4}}{\longleftarrow} & \cdots\\
& \nwarrow\quad\swarrow &  & \nwarrow\quad\swarrow &  & \nwarrow\quad\swarrow
&  &  & \\
&  im(\lambda_{1}) & \twoheadleftarrow &  im(\lambda_{2}) & \twoheadleftarrow
&  im(\lambda_{3}) & \twoheadleftarrow & \cdots &
\end{array}
\]
The $im(\lambda_{i})$'s are clearly finitely generated but may not be finitely
presented. We will use this diagram to produce a new sequence with the desired properties.

For each $i\geq1$, let $\left\{  g_{j}^{i}\right\}  _{j=1}^{n_{i}}$ be a
generating set for $G_{i}$, and choose $\left\{  h_{j}^{i}\right\}
_{j=1}^{n_{i}}$ $\subset im(\lambda_{i+1})$ so that $\lambda_{i}(g_{j}%
^{i})=\lambda_{i}(h_{j}^{i})$. 

\textbf{Note}\qua The superscripts are indices,
not powers.

Let $H_{i}$ $\vartriangleleft G_{i}$ be the normal closure of the
set $S_{i}=\left\{  g_{j}^{i}\left(  h_{j}^{i}\right)  ^{-1}\right\}
_{j=1}^{n_{i}}$, define $G_{i}^{\prime}$ to be $G_{i}/H_{i}$, and let
$q_{i}\co G_{i}\rightarrow G_{i}^{\prime}$ be the quotient map. Since
$S_{i+1}\subset\ker(\lambda_{i+1})$ we get induced homomorphisms
$\lambda_{i+1}^{\prime}\co G_{i+1}^{\prime}\rightarrow G_{i}$ . Define $\mu
_{i+1}=q_{i}\circ\lambda_{i+1}^{\prime}\co G_{i+1}^{\prime}\rightarrow
G_{i}^{\prime}$ to obtain the commuting diagram:
\[%
\begin{array}
[c]{cccccccc}%
& G_{1} & \overset{\lambda_{2}}{\longleftarrow} & G_{2} & \overset{\lambda
_{3}}{\longleftarrow} & G_{3} & \overset{\lambda_{4}}{\longleftarrow} &
\cdots\\
\qquad\overset{q_{1}}{\swarrow} &  & \overset{\lambda_{2}^{\prime}}{\nwarrow
}\quad\overset{q_{2}}{\swarrow} &  & \overset{\lambda_{3}^{\prime}}{\nwarrow
}\quad\overset{q_{3}}{\swarrow} &  &  & \\
G_{1}^{\prime} & \overset{\mu_{2}}{\twoheadleftarrow} & G_{2}^{\prime} &
\overset{\mu_{3}}{\twoheadleftarrow} & G_{3}^{\prime} & \overset{\mu_{4}%
}{\twoheadleftarrow} & \cdots &
\end{array}
\]
Since each $G_{i}^{\prime}$ is generated by $\left\{  q_{i}(g_{j}%
^{i})\right\}  _{j=1}^{n_{i}}$ and each $q_{i}(g_{j}^{i})$ has a preimage
$h_{j}^{i}$ $\in G_{i+1}$ under the map $q_{i}\circ\lambda_{i+1}$, it follows
(from the commutativity of the diagram) that each $\mu_{i+1}$ is surjective.
Lastly, each $G_{i}^{\prime}$ has a finite presentation which may be obtained
from a finite presentation for $G_{i}$ by adding relators corresponding to the
elements of $S_{i}$.
\end{proof}

\section{Ends of manifolds: definitions and background
information\label{background}}

In this section we review some standard notions involved in the study of
non-compact manifolds and complexes. Since the terminology and notation used
in this area are by no means standardized, the reader should be careful when
consulting other sources. The remarks at the end of the section addresses a
portion of this issue.

The symbol $\approx$ will denote homeomorphisms; $\simeq$ will denote
homotopic maps or homotopy equivalent spaces. A manifold $M^{n}$ is
\emph{open} if it is non-compact and has no boundary. We say that $M^{n}$ is
\emph{one ended} if complements of compacta in $M^{n}$ contain exactly one
unbounded component. For convenience, we restrict our attention to one ended
manifolds. In addition, we will work in the PL category. Equivalent results in
the smooth and topological categories may be obtained in the usual ways.
Results may be generalized to spaces with finitely many ends by considering
one end at a time.

A set $U\subset M^{n}$ is a \emph{neighborhood of infinity} if $\overline
{M^{n}-U}$ is compact; $U$ is a \emph{clean neighborhood of infinity }if it is
also a PL submanifold with bicollared boundary. It is easy to see that each
neighborhood $U$ of infinity contains a clean neighborhood $V$ of
infinity---just let $V=M^{n}-$ $\overset{\circ}{N}$ where $N$ is a regular
neighborhood of a polyhedron containing $\overline{M^{n}-U}$. We may also
arrange that $V$ be connected by discarding all of its compact components.
Thus we have:

\begin{lemma}
Each one ended open manifold $M^{n}$ contains a sequence $\left\{
U_{i}\right\}  _{i=0}^{\infty}$ of clean connected neighborhoods of infinity
with $U_{i+1}\subset\overset{\circ}{U}_{i}$ for all $i\geq0$, and
$\bigcap_{i=0}^{\infty}U_{i}=\emptyset$.
\end{lemma}

\noindent A sequence of the above type will be called \emph{neat}. In the
future, all neighborhoods of infinity are assumed to be clean and connected
and sequences of these neighborhoods are neat.

We say that $M^{n\text{ }}$is \emph{inward tame }at infinity if, for
arbitrarily small neighborhoods of infinity $U$, there exist homotopies
$H\co U\times\left[  0,1\right]  \rightarrow U$ such that $H_{0}=id_{U}$ and
$\overline{H_{1}\left(  U\right)  }$ is compact. Thus inward tameness means
that neighborhoods of infinity can be pulled into compact subsets of themselves.

Recall that a CW complex $X$ is \emph{finitely dominated} if there exists a
finite complex $K$ and maps $u\co X\rightarrow K$ and $d\co K\rightarrow X$ such
that $d\circ u\simeq id_{X}$. It is easy to see that $X$ \ is finitely
dominated if and only if it may be homotoped into a compact subset of itself.
Hence, our manifold $M^{n}$ is inward tame if and only if arbitrarily small
neighborhoods of infinity are finitely dominated. This characterization of
``inward tameness'' will be useful to us later.

Next we study the \emph{fundamental group system at the end} of $M^{n}$. Begin
with a neat sequence $\left\{  U_{i}\right\}  _{i=0}^{\infty}$ of
neighborhoods of infinity and basepoints $p_{i}\in U_{i}$. For each $i\geq1$,
choose a path $\alpha_{i}\subset U_{i-1}$ connecting $p_{i}$ to $p_{i-1}$.
Then, for each $i\geq1$, let $\lambda_{i}\co \pi_{1}\left(  U_{i},p_{i}\right)
\rightarrow\pi_{1}\left(  U_{i-1},p_{i-1}\right)  $ be the homomorphism
induced by inclusion followed by the change of basepoint isomorphism
determined by $\alpha_{i}$. Suppressing basepoints, this gives us an inverse
sequence:
\[
\pi_{1}\left(  U_{0}\right)  \overset{\lambda_{1}}{\longleftarrow}\pi
_{1}\left(  U_{1}\right)  \overset{\lambda_{2}}{\longleftarrow}\pi_{1}\left(
U_{2}\right)  \overset{\lambda_{3}}{\longleftarrow}\pi_{1}\left(
U_{3}\right)  \overset{\lambda_{3}}{\longleftarrow}\cdots
\]
Provided this sequence is semistable, one can show that its pro-equivalence
class does not depend on any of the choices made above. We then denote the
pro-equivalence class of this sequence by $\pi_{1}\left(  \varepsilon\left(
M^{n}\right)  \right)  $. (In the absence of semistability, the choices become
part of the data.) We will denote the inverse limit of the above sequence by
$\pi_{1}\left(  \infty\right)  $.\smallskip

\noindent\textbf{Note}\qua For our purposes, $\pi_{1}\left(  \infty\right)  $
will only be used as a (rather trivial) convenience when $\pi_{1}\left(
\varepsilon\left(  M^{n}\right)  \right)  $ is stable. Otherwise, we work with
the inverse sequence.$\smallskip$

Before moving to a new topic, notice that the same procedure may be used to
define $\pi_{k}\left(  \varepsilon\left(  M^{n}\right)  \right)  $ for $k>1$.

We now understand the first two conditions in Theorem \ref{sieb}, and begin to
look at the third. If $\Lambda$ is a ring, we say that two finitely generated
projective $\Lambda$--modules $P$ and $Q$ are \emph{stably equivalent} if there
exist finitely generated free $\Lambda$--modules $F_{1\text{ }}$ and $F_{2}$
such that $P\oplus F_{1}\cong Q\oplus F_{2}$. The stable equivalence classes
of finitely generated projective modules form a group $\widetilde{K}%
_{0}\left(  \Lambda\right)  $ under direct sum. Then $P$ represents the
trivial element of $\widetilde{K}_{0}\left(  \Lambda\right)  $ if and only if
it is \emph{stably free}, ie, there exists a finitely generated free
$\Lambda$--module $F$ such that $P\oplus F$ is free. In \cite{Wa1}, Wall shows
that each finitely dominated $X$ determines a well-defined element
$\sigma\left(  X\right)  \in\widetilde{K}_{0}\left(  \mathbb{Z}\left[  \pi
_{1}X\right]  \right)  $ which vanishes if and only if $X$ has the homotopy
type of a finite complex. When an open one ended manifold $M^{n}$ ($n\geq6$)
satisfies Conditions 1 and 2 of Theorem \ref{sieb}, Siebenmann isolated a
single obstruction (which we have denoted $\sigma_{\infty}\left(
M^{n}\right)  $) in $\widetilde{K}_{0}\left(  \mathbb{Z}\left[  \pi_{1}%
(\infty)\right]  \right)  $ to finding an open collar neighborhood of
infinity. In addition he observed that, up to sign, his obstruction is just
the Wall obstruction of an appropriately chosen neighborhood of infinity. One
upshot of this observation (requiring use of Siebenmann's Sum Theorem for the
Finiteness Obstruction---see \cite{Si1} or \cite{Fe2}) is that $\sigma
_{\infty}\left(  M^{n}\right)  $ vanishes if and only if all clean
neighborhoods of infinity in $M^{n}$ have finite homotopy types.

When $\pi_{1}\left(  \varepsilon\left(  M^{n}\right)  \right)  $ is not
stable, the definition of $\sigma_{\infty}\left(  M^{n}\right)  $ becomes
somewhat more complicated. Instead of measuring the obstruction in a single
neighborhood of infinity, it will lie in the group $\widetilde{K}_{0}\left(
\pi_{1}\left(  \varepsilon\left(  M^{n}\right)  \right)  \right)
\equiv{\displaystyle\lim_{\textstyle\longleftarrow}}\widetilde{K}_{0}\left(  \mathbb{Z}\left[  \pi
_{1}U_{i}\right]  \right)  $, where $\left\{  U_{i}\right\}  $ is a neat
sequence of neighborhoods of infinity. Then $\sigma_{\infty}(M^{n})$ may be
identified with the element $\left(  -1\right)  ^{n}(\sigma(U_{0}%
),\sigma(U_{1}),\sigma(U_{2}),\cdots)$, with $\sigma\left(  U_{i}\right)  $
being the Wall finiteness obstruction for $U_{i}$. Again, this obstruction
vanishes if and only if all clean neighborhoods of infinity in $M^{n}$ have
finite homotopy types. When $\pi_{1}\left(  \varepsilon(M^{n}\right)  )$ is
stable, this definition of $\sigma_{\infty}(M^{n})$ reduces to the one
discussed above. When $n\geq6$ and $\pi_{1}\left(  \varepsilon\left(
M^{n}\right)  \right)  $ is semistable, we will see $\sigma_{\infty}(M^{n})$
arise naturally---without reference to the Wall finiteness obstruction---as an
obstruction to pseudo-collar\-abil\-ity (see Section \ref{n-2}). For a more
general treatment of this obstruction---which, among other things, shows that
$\widetilde{K}_{0}\left(  \pi_{1}\left(  \varepsilon\left(  M^{n}\right)
\right)  \right)  $ and $\sigma_{\infty}(M^{n})$ are independent of the choice
of $\left\{  U_{i}\right\}  $---we refer the reader to \cite{CS}.

\begin{remark}
Our use of the phrase ``inward tame'' is not standard. In \cite{CS} the same
notion is simply called ``tame'', while in \cite{Si1}, ``tame'' means ``inward
tame and $\pi_{1}$--stable''. Quinn and others (see, for example, \cite{HR})
have given ``tame'' a different and inequivalent meaning which involves
pushing neighborhoods of infinity toward the end of the space, while referring
to our brand of tameness as ``reverse tameness''. We hope that by referring to
our version of tameness as ``inward tame'' and Quinn's version as ``outward
tame'' we can avoid some confusion.
\end{remark}

\begin{remark}
One should be careful not to interpret the symbol $\sigma_{\infty}\left(
M^{n}\right)  $ as the Wall finiteness obstruction $\sigma\left(
M^{n}\right)  $ of the manifold $M^{n}$. Indeed, $M^{n}$ can have finite
homotopy type even when its neighborhoods of infinity do not. (The Whitehead
contractible $3$--manifold is one well-known example.) This situation can arise
even when $\pi_{1}\left(  \varepsilon\left(  M^{n}\right)  \right)  $ is stable.
\end{remark}

\section{Pseudo-collars and the Main Theorem\label{pseudo}}

Recall that a manifold $U^{n}$ with compact boundary is an \emph{open collar}
if $U^{n}\approx\partial U^{n}\times\lbrack0,\infty)$; it is a \emph{homotopy
collar} if the inclusion $\partial U^{n}\hookrightarrow U^{n}$ is a homotopy
equivalence. If $U^{n}$ is a homotopy collar which contains arbitrarily small
homotopy collar neighborhoods of infinity, then we call $U^{n}$ a
\emph{pseudo-collar}. We say that an open $n$--manifold $M^{n}$ is
\emph{collarable} if it contains an open collar neighborhood of infinity, and
that $M^{n}$ is \emph{pseudo-collarable} if it contains a pseudo-collar
neighborhood of infinity. The following easy example is useful to keep in mind.

\begin{example}
Let $M^{n}$ be a contractible $n$--manifold and $B^{n}\subset M^{n}$ a
standardly embedded $n$--ball. Then $U=M^{n}-\overset{\circ}{B^{n}}$ is a
homotopy collar; however, in general $M^{n}$ need not be pseudo-collarable
(see Example \ref{whitehead}).
\end{example}

\begin{remark}
A standard duality argument guarantees that any connected homotopy collar
(hence any connected pseudo-collar) is one ended. See, for example, \cite{Si2}.
\end{remark}

When discussing collars, some complementary notions are useful. A compact
codimension $0$ submanifold $C$ of an open manifold $M^{n}$ is called a
\emph{core} if $C\hookrightarrow M^{n}$ is a homotopy equivalence; it is
called a \emph{geometric core} if $M^{n}-\overset{\circ}{C}$ is a homotopy
collar; and it is called an \emph{absolute core }if $M^{n}-\overset{\circ}{C}$
is an open collar. The following is immediate.

\begin{proposition}
Let $M^{n}$ be a one ended open $n$--manifold. Then:

\begin{enumerate}
\item $M^{n}$ is collarable if and only if $M^{n}$ contains an absolute core
(hence, arbitrarily large absolute cores), and

\item $M^{n}$ is pseudo-collarable if and only if $M^{n}$ contains arbitrarily
large geometric cores.
\end{enumerate}
\end{proposition}

\begin{example}
\label{whitehead}The Whitehead contractible $3$--manifold $M^{3}$ is not
pseudo-coll\-arable. Indeed, if $M^{3}$ were pseudo-collarable, it would contain
arbitrarily large geometric cores each of which---by standard $3$--manifold
topology---would be a $3$--ball. But then $M^{3}$ would be a monotone union of
open $3$--balls, and hence, homeomorphic to $\mathbb{R}^{3}$ by \cite{Br} or by
an application of the Combinatorial Annulus Theorem (Corollary 3.19 of
\cite{RS}).
\end{example}

On the positive side we have:

\begin{example}
Although they are not collarable, the exotic universal covering spaces
produced by Davis in \cite{Da} are pseudo-collarable. If $M^{n}$ is one of
these covering spaces with compact contractible manifold $C^{n}$ as a
``fundamental chamber'', then $M^{n}$ contains arbitrarily large geometric
cores homeomorphic to finite sums
\[
C^{n}\#_{\partial}C^{n}\#_{\partial}\cdots\#_{\partial}C^{n}%
\]
(with increasing numbers of summands). Here $\#_{\partial}$ denotes a
``boundary connected sum'', ie, the union of two $n$--manifolds with boundary
along boundary $\left(  n-1\right)  $--disks. However, $M^{n}$ is not
collarable since $\pi_{1}\left(  \varepsilon\left(  M^{n}\right)  \right)  $
is not stable---in fact $\pi_{1}\left(  \varepsilon\left(  M^{n}\right)
\right)  $ may be represented by the sequence
\[
G\leftarrow G\ast G\leftarrow(G\ast G)\ast G\leftarrow\left(  G\ast G\ast
G\right)  \ast G\leftarrow\cdots
\]
where $G=\pi_{1}\left(  \partial C^{n}\right)  $ and each homomorphism is
projection onto the first term. It is interesting to note that this sequence
is perfectly semistable.
\end{example}

A compact cobordism $\left(  W^{n},M^{n-1},N^{n-1}\right)  $ is a
\emph{one-sided h-cobordism} if one (but not necessarily both) of the
inclusions $M^{n-1}\hookrightarrow W^{n}$ or $N^{n-1}\hookrightarrow W^{n}$ is
a homotopy equivalence. The following property of one-sided $h$-cobordisms is
a well-known consequence of duality (see, for example, Lemma 2.5 of \cite{DT1}).

\begin{lemma}
\label{perfect}Let $\left(  W^{n},M^{n-1},N^{n-1}\right)  $ be a compact
connected one-sided $h$-cob\-ordism with $M^{n-1}\overset{\simeq}%
{\hookrightarrow}W^{n}$. Then the inclusion induced homomorphism\break $\pi
_{1}\left(  N^{n-1}\right)  \rightarrow\pi_{1}\left(  W^{n}\right)  $ is
surjective and has perfect kernel.
\end{lemma}

Non--trivial one-sided $h$-cobordisms are plentiful. In fact, if we
are given a closed $\left( n-1\right) $--manifold $N^{n-1}$
($n\geq6$), a finitely presented group $G$, and a homomorphism $\mu\co
\pi_{1}\left( N^{n-1}\right) \twoheadrightarrow G$ with perfect
kernel, then the ``Quillen plus construction'' (see \cite{Qu} or Sections
11.1 and 11.2 of \cite{FQ}) produces a one-sided $h$-cobordism $\left(
W^{n},M^{n-1},N^{n-1}\right) $ with $\pi_{1}\left( W^{n}\right) \cong
G$ and $$\ker\left( \pi_{1}\left( N^{n-1}\right)
\rightarrow\pi_{1}\left( W^{n}\right) \right) =\ker\left( \mu\right)
.$$
The role played by one-sided $h$-cobordisms in the study of pseudo-collars is
clearly illustrated by the following easy proposition.

\begin{proposition}
\label{decomposition}Let $\left\{  \left(  W_{i},M_{i},N_{i}\right)  \right\}
_{i=1}^{\infty}$ be a collection of one-sided h-cobord\-isms with $M_{i}%
\overset{\simeq}{\hookrightarrow}W_{i}$, and suppose that for each $i\geq1 $
there is a homeomorphism $h_{i}\co N_{i}\rightarrow M_{i+1}$. Then the adjunction
space
\[
U=W_{1}\cup_{h_{1}}W_{2}\cup_{h_{2}}W_{3}\cup_{h_{3}}\cdots
\]
is a pseudo-collar. Conversely, every pseudo-collar may be expressed as a
countable union of one-sided h-cobordisms in this manner.
\end{proposition}

\begin{proof}
For the forward implication, we begin by observing that $U$ is a homotopy
collar. First note that $\partial U=M_{1}\hookrightarrow W_{1}\cup_{h_{1}%
}\cdots\cup_{h_{1}}W_{k}$ is a homotopy equivalence for any finite $k$. A
direct limit argument then shows that $\partial U\hookrightarrow U$ is a
homotopy equivalence. Alternatively, we may observe that $\pi_{\ast}\left(
U,\partial U\right)  \equiv0$ and apply the Whitehead theorem. To see that $U$
is a pseudo-collar we apply the same argument to the subsets $U_{i}%
=W_{i+1}\cup_{h_{i+1}}W_{i+2}\cup_{h_{i+2}}W_{i+3}\cup_{h_{i+3}}\cdots$.

For the converse, assume that $U$ is a pseudo-collar. Choose a homotopy collar
$U_{1}\subset\overset{\circ}{U}$ and let $W_{1}=U-\overset{\circ}{U}_{1}$.
Then $\partial U\overset{\simeq}{\hookrightarrow}W_{1}$, so $\left(
W_{1},\partial U,\partial U_{1}\right)  $ is a one-sided $h$-cobordism. Next
choose a homotopy collar $U_{2}$ $\subset\overset{\circ}{U}_{1}$ and let
$W_{2}=U_{1}-\overset{\circ}{U}_{2}$. Repeating this procedure gives the
desired result. See Figure 1.%
\begin{figure}
[ht!]
\begin{center}
\includegraphics[width=4.in]{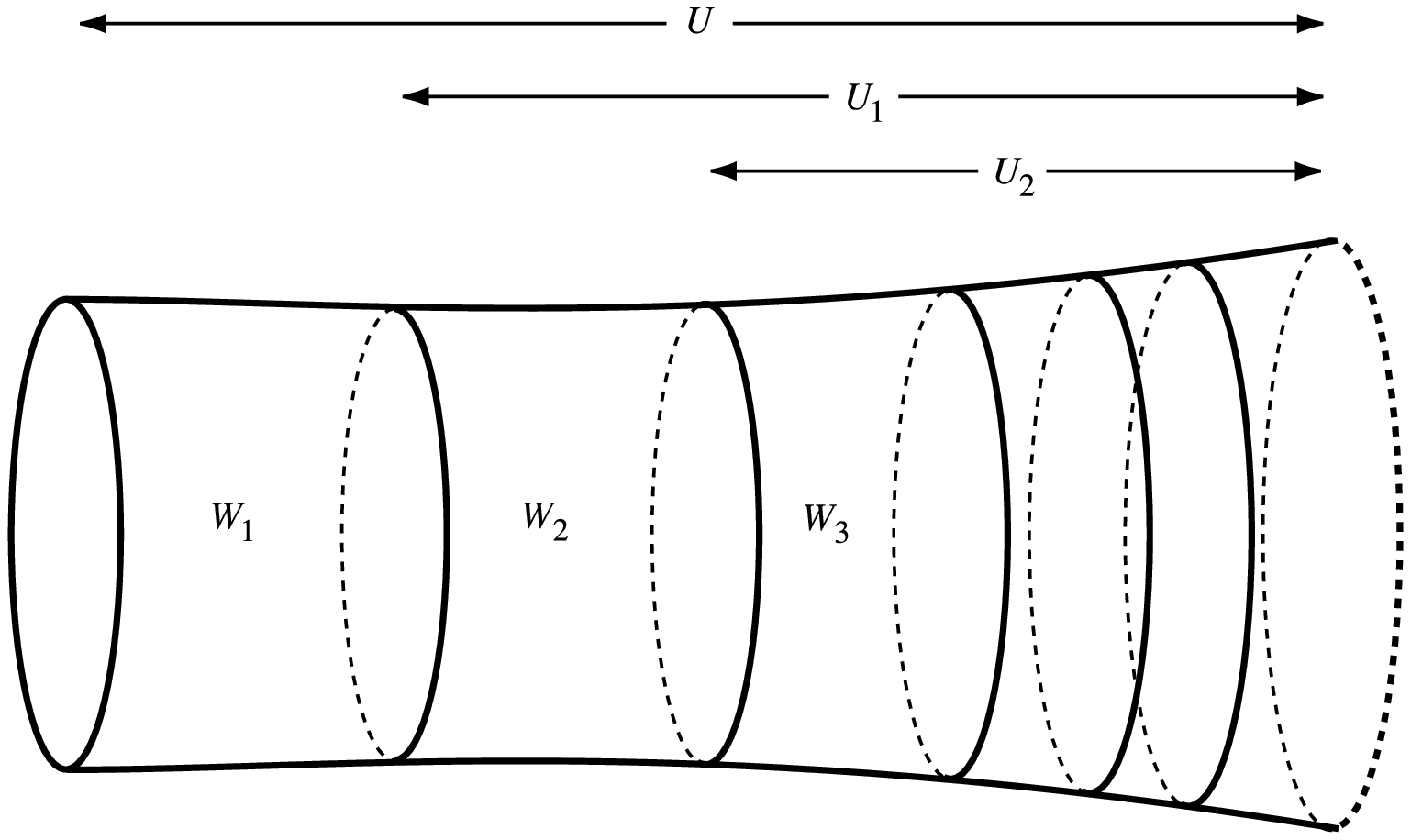}
\caption{}
\end{center}
\end{figure}
\end{proof}

The next result provides a striking similarity between pseudo-collars and
genuine open collars. It follows immediately from Proposition
\ref{decomposition} and the main result of \cite{DT2} which shows that
one-sided $h$-cobordisms in dimensions $\geq6$ may be ``laminated''.

\begin{proposition}
\label{laminated}Let $U^{n}$ be a pseudo-collar ($n\geq6$). Then there exists
a proper continuous surjection $p\co U^{n}\rightarrow\lbrack0,\infty)$ with the
following properties.

\begin{enumerate}
\item $p^{-1}\left(  0\right)  =\partial U^{n},$

\item  each $p^{-1}\left(  r\right)  $ is a closed $\left(  n-1\right)
$--manifold with the same $\mathbb{Z}$--homology as $\partial U^{n}$, and

\item $p^{-1}\left(  r\right)  $ is nicely embedded, ie, has a product
neighborhood in $U^{n}$, for $r\neq1,2,3,\cdots$.
\end{enumerate}
\end{proposition}

Our next result provides the fundamental conditions necessary for
pseudo-\allowbreak collar\-abil\-\allowbreak ity.

\begin{theorem}
\label{necessary}Suppose a one ended open manifold $M^{n}$ is
pseudo-collarable. Then

\begin{enumerate}
\item $M^{n}$ is inward tame at infinity,

\item $\pi_{1}(\varepsilon(M^{n}))$ is perfectly semistable, and

\item $\sigma_{\infty}\left(  M^{n}\right)  =0\in\widetilde{K}_{0}\left(
\pi_{1}\left(  \varepsilon(M^{n}\right)  \right)  )$.
\end{enumerate}
\end{theorem}

\begin{proof}
Properties 1 and 3 follow easily from the definition of pseudo-collar; while
Property 2 is obtained from Proposition \ref{decomposition} and Lemma
\ref{perfect}.
\end{proof}

One might hope that the above conditions are also sufficient for $M^{n}$
($n\geq6$) to be pseudo-collarable. This would be an ideal generalization of
Theorem \ref{sieb}; but, although we have not ruled it out, we are thus far
unable to prove it. Our main result---which, for easy reference, we now
restate---requires an additional hypothesis and one additional dimension.

\begin{theorem}
[Main Existence Theorem]\label{met}A one ended open $n$--manifold $M^{n}$
($n\geq7$) is pseudo-collarable provided each of the following is satisfied:

\begin{enumerate}
\item $M^{n}$ is inward tame at infinity,

\item $\pi_{1}(\varepsilon(M^{n}))$ is perfectly semistable,

\item $\sigma_{\infty}\left(  M^{n}\right)  =0\in\widetilde{K}_{0}\left(
\pi_{1}\left(  \varepsilon(M^{n}\right)  \right)  )$, and

\item $\pi_{2}(\varepsilon(M^{n}))$ is semistable.
\end{enumerate}
\end{theorem}

\begin{remark}
Several interesting classes of manifolds are known to satisfy some or all of
the conditions\ in the above theorems, thus making them ideal candidates for
pseudo-collar\-abil\-ity. We mention a few of them. 

{(a)}\qua We already
know that the exotic universal coverings of \cite{Da} are pseudo-collarable,
and therefore satisfy Conditions 1--3. It can also be shown that they satisfy
Condition 4. 

{(b)}\qua Every piecewise flat $CAT\left(  0\right)  $
manifold satisfies Conditions 1--3. Some of the most interesting of these---the
exotic universal covers produced by Davis and Januszkiewicz in \cite{DJ}%
---also satisfy Condition 4 (and are therefore pseudo-collarable). 

{(c)}\qua A more general class of open $n$--manifolds which are of current interest are
those admitting $\mathcal{Z}$--compactifications (see \cite{AG}, \cite{Be},
\cite{FW} and \cite{CP} for discussions).\ These manifolds satisfy Conditions
1 and 3, and also have semistable fundamental groups at infinity (whether
these are perfectly semistable is unknown).
\end{remark}

Most of the remainder of this paper is devoted to proving the Main Existence Theorem.

\section{Proof of the Main Existence Theorem: an outline\label{outline}}

Let $M^{n}$ be a 1--ended open manifold and $U$ a connected clean neighborhood
of infinity. According to \cite{Si1}, $U$ is a $0$\emph{--neighborhood of
infinity} if $\partial U$ is connected. Under the assumption that $\pi
_{1}\left(  \varepsilon\left(  M^{n}\right)  \right)  $ is stable, \cite{Si1}
then defines $U$ to be a \emph{1--neighborhood of infinity }provided it is a
$0$--neighborhood infinity and both $\pi_{1}\left(  \infty\right)
\rightarrow\pi_{1}\left(  U\right)  $ and $\pi_{1}\left(  \partial U\right)
\rightarrow\pi_{1}\left(  U\right)  $ are isomorphisms. For $k\geq2 $, $U$ is
a \emph{k--neighborhood of infinity }if it is a $1$--neighborhood of infinity
and $\pi_{i}(U,\partial U)=0$ for $i\leq k$.

We may now describe Siebenmann's proof of Theorem \ref{sieb}. Beginning with a
neat sequence $\left\{  U_{i}\right\}  $ of neighborhoods of infinity, perform
geometric alterations to obtain a neat sequence of $0$--neighborhoods of
infinity. This is easy---given a $U_{i}$ with non-connected boundary, choose
finitely many disjoint properly embedded arcs in $U_{i}$ connecting the
components of $\partial U_{i}$. Then ``drill out'' regular neighborhoods of
these arcs to connect up the boundary components, thus obtaining a
$0$--neighborhood $U_{i}^{\prime}\subset U_{i}$. After passing to a subsequence
(if necessary) to maintain the ``nestedness'' condition, we have the desired
sequence. Assuming then that $\left\{  U_{i}\right\}  $ is a neat sequence of
$0$--neighborhoods of infinity and that Conditions 1 and 2 of Theorem
\ref{sieb} are both satisfied, convert the $U_{i}$'s into $1$--neighborhoods of
infinity. This stage of the proof is more complicated. We view it as the first
of three major steps in obtaining Theorem \ref{sieb}. Some algebra (Lemmas
\ref{gp3} and \ref{gp4}) is required, neighborhoods of arcs are drilled out,
and neighborhoods of disks are ``traded''---sometimes removed and sometimes
added. Ultimately one obtains a neat sequence of $1$--neighborhoods of
infinity. Next, in the middle step of the proof, the $U_{i} $'s are
inductively improved until they are $\left(  n-3\right)  $--neighborhoods. The
key tools here are: general position, handle theory, and Lemma \ref{fg}. The
final step in Siebenmann's proof is to improve $\left(  n-3\right)
$--neighborhoods of infinity to $\left(  n-2\right)  $--neighborhoods---which
turn out to be open collars. This step is very delicate. More algebra is
required, the need for Condition 3 becomes clear, and $\pi_{1}$--stability
plays a crucial role.

To a large extent, the proof of our Main Existence Theorem is a careful
reworking of \cite{Si1}. In fact, the reader will find Siebenmann's proof
properly embedded in ours. However, since the $\pi_{1}$--stability hypothesis
so thoroughly permeates \cite{Si1}, a great deal of revision and
generalization is necessary. First, we define a \emph{generalized
1--neighborhood of infinity }to be a $0$--neighborhood of infinity $U$ with the
property that $\pi_{1}\left(  \partial U\right)  \rightarrow\pi_{1}\left(
U\right)  $ is an isomorphism. Then for $k\geq2$, a \emph{generalized
k--neighborhood of infinity }is a generalized $1$--neighborhood of infinity with
the property that $\pi_{i}(U,\partial U)=0$ for $i\leq k$. The point here is
that, when $\pi_{1}\left(  \varepsilon\left(  M^{n}\right)  \right)  $ is not
stable, there is no ``preferred fundamental group'' for our neighborhoods of
infinity. Later we will see that there are sometimes ``preferred sequences of
fundamental groups''. To avoid confusion, we will often refer to the
$k$--neighborhoods of infinity defined earlier as \emph{strong k--neighborhoods
of infinity}. Of course, this only makes sense when $\pi_{1}\left(
\varepsilon\left(  M^{n}\right)  \right)  $ is stable.

We break our proof into the same three major steps as above. In the first step
(Section \ref{1}) we obtain neat sequences of generalized $1$--neighborhoods of
infinity. For this, only Condition 1 of the Main Existence Theorem is
required; however, given additional assumptions about $\pi_{1}\left(
\varepsilon\left(  M^{n}\right)  \right)  $ (eg, stability, semistability,
or perfect semistability), we show how these may be incorporated. The middle
step of the proof (Section \ref{n-3}) requires the least revision of
\cite{Si1}. Only Condition 1 is needed to obtain a neat sequence of
generalized $\left(  n-3\right)  $--neighborhoods of infinity. As before,
additional assumptions on the fundamental group at infinity can be
incorporated into this step. The final step (Section \ref{n-2}) is the most
difficult. In order to make any progress beyond generalized $\left(
n-3\right)  $--neighborhoods of infinity, it becomes necessary to assume that
$\pi_{1}\left(  \varepsilon\left(  M^{n}\right)  \right)  $ is semistable (a
part of Condition 2). We show that a neat sequence of generalized $\left(
n-2\right)  $--neighborhoods, with $\pi_{1}$--semistability appropriately built
in, determines a pseudo-collar structure; hence, obtaining generalized
$\left(  n-2\right)  $--neighborhoods is our goal. In our attempt to mimic
Siebenmann, we rediscover the $\widetilde{K}_{0}$--obstruction much as it
appeared in \cite{Si1}. The difference is that, since $\pi_{1}\left(
U_{i}\right)  $ changes with $i$, so must the $\widetilde{K}_{0}$--obstruction.
Hence, our obstruction becomes a \emph{sequence} of obstructions. When this
obstruction dies, most of the algebraic and handle theoretic steps from
\cite{Si1} may be duplicated. Unfortunately, at the last instant---a final
application of the Whitney Lemma---the lack of $\pi_{1}$--stability creates
major problems.\ To complete the proof in the non-stable situation, we are
forced to develop a new strategy. It is only here that we require Condition 4
and the ``perfect'' part of Condition 2.

\section{Obtaining generalized $1$--neighborhoods of infinity\label{1}}

In this section we show how to obtain a neat sequence $\left\{  U_{i}\right\}
$ of generalized $1$--neighborhoods of infinity in a one ended open
$n$--manifold when $n\geq5.$ This requires only that $M^{n}$ be inward tame at
infinity. (In fact, it would be enough to assume that clean neighborhoods of
infinity have finitely presentable fundamental groups.) In addition we show
that, when $\pi_{1}(\varepsilon(M^{n}))$ is pro-equivalent to certain
preferred inverse sequences of surjections, we can make our sequence $\left\{
\pi_{1}(U_{i})\right\}  $ isomorphic to corresponding subsequences \ This
covers situations where $\pi_{1}(\varepsilon(M^{n}))$ is stable, semistable
and perfectly semistable.

\begin{lemma}
\label{lemma1ngh}Let $M^{n}$ ($n\geq5$) be a one ended open $n$--manifold which
is inward tame at infinity and let $V$ be a $0$--neighborhood of infinity. Then
$V$ contains a generalized $1$--neighborhood $U$ of infinity with the property
that $\pi_{1}\left(  U\right)  \rightarrow\pi_{1}\left(  V\right)  $ is an isomorphism.
\end{lemma}

\begin{proof}
First we construct a $0$--neighborhood $V^{\prime}\subset V$ so that $\pi
_{1}\left(  \partial V^{\prime}\right)  \rightarrow\pi_{1}\left(  V^{\prime
}\right)  $ is surjective and $\pi_{1}\left(  V^{\prime}\right)
\rightarrow\pi_{1}\left(  V\right)  $ is an isomorphism.

Since $V$ is finitely dominated, $\pi_{1}\left(  V\right)  $ is finitely
generated, so we may choose a finite collection $\left\{  \alpha_{1}%
,\alpha_{2},\cdots,\alpha_{k}\right\}  $ of disjoint properly embedded p.l.
arcs in $V$ so that $\pi_{1}(\partial V\cup(\bigcup_{i=1}^{k}\alpha
_{i}))\rightarrow\pi_{1}(V)$ is surjective. Choose a collection $\left\{
N_{i}\right\}  _{i=1}^{k}$ of disjoint regular neighborhoods of the
$\alpha_{i}$'s in $V$ and let
\[
V^{\prime}=\overline{V-\bigcup\nolimits_{i=1}^{k}N_{i}}.
\]
Clearly $\pi_{1}\left(  \partial V^{\prime}\right)  $ (and thus $\pi
_{1}\left(  V^{\prime}\right)  $) surjects onto $\pi_{1}\left(  V\right)  $;
moreover, since disks in $V$ may be pushed off the $N_{i}$'s, then $\pi
_{1}\left(  V^{\prime}\right)  \rightarrow\pi_{1}\left(  V\right)  $ is also injective.

Next we modify $V^{\prime}$ to be a generalized $1$--neighborhood. Since
$V^{\prime}$ is finitely dominated, Lemma \ref{gp3} implies that $\pi
_{1}\left(  V^{\prime}\right)  $ is finitely presentable. Hence, by Lemma
\ref{gp4}, $\ker(\pi_{1}\left(  \partial V^{\prime}\right)  \rightarrow\pi
_{1}\left(  V^{\prime}\right)  )$ is the normal closure of a finite set of
elements. Let $\left\{  \beta_{1},\beta_{2},\cdots,\beta_{r}\right\}  $ be a
collection of pairwise disjoint embedded loops in $\partial V^{\prime}$
representing these elements, then choose $\left\{  D_{1},D_{2},\cdots
,D_{r}\right\}  $ a pairwise disjoint collection of properly embedded 2--disks
in $V^{\prime}$ with $\partial D_{i}=\beta_{i}$ for each $i$. Let $\left\{
P_{1},P_{2},\cdots,P_{r}\right\}  $ be a pairwise disjoint collection of
regular neighborhoods of the $D_{i}$'s in $V^{\prime}$ and define
\[
U=\overline{V^{\prime}-\bigcup\nolimits_{i=1}^{r}P_{i}}.
\]
By VanKampen's theorem $\pi_{1}\left(  \partial V^{\prime}\cup(\bigcup
_{i=1}^{r}P_{i})\right)  \rightarrow\pi_{1}\left(  V^{\prime}\right)  $ is an
isomorphism, and by general position $\pi_{1}\left(  \partial U\right)
\rightarrow\pi_{1}\left(  \partial V^{\prime}\cup(\bigcup_{i=1}^{r}%
P_{i})\right)  $ and $\pi_{1}\left(  U\right)  \rightarrow\pi_{1}\left(
V^{\prime}\right)  $ are isomorphisms. It follows that $U$ is a generalized
$1$-- neighborhood of infinity and $\pi_{1}\left(  U\right)  \rightarrow\pi
_{1}\left(  V\right)  $ is an isomorphism.
\end{proof}

Combining the above lemma with the method described in the previous section
for obtaining $0$--neighborhoods of infinity gives:

\begin{corollary}
\label{cor1ngh}Every one ended open $n$--manifold ($n\geq5$) that is inward
tame at infinity contains a neat sequence of generalized $1$--neighborhoods of infinity.
\end{corollary}

\begin{lemma}
\label{lemma1seq}Let $M^{n}$ ($n\geq5$) be a one ended $n$--manifold that is
inward tame at infinity and suppose the fundamental group system $\pi
_{1}(\varepsilon(M^{n}))$ is pro-equivalent to an inverse sequence
$\mathcal{G}\co $ $G_{1}\twoheadleftarrow G_{2}\twoheadleftarrow G_{3}%
\twoheadleftarrow\cdots$ of finitely presentable groups and surjections. Then
there is a neat sequence $\left\{  U_{i}\right\}  _{i=1}^{\infty}$ of
$1$--neighborhoods of infinity so that the inverse sequence $\pi_{1}%
(U_{1})\leftarrow\pi_{1}(U_{2})\leftarrow\pi_{1}(U_{3})\leftarrow\cdots$ is
isomorphic to a subsequence of $\mathcal{G}$.
\end{lemma}

\begin{proof}
By the hypothesis and Corollary \ref{cor1ngh}, there exists a neat sequence
$\left\{  V_{i}\right\}  $ of generalized $1$--neighborhoods of infinity, a
subsequence $G_{k_{1}}\twoheadleftarrow G_{k_{2}}\twoheadleftarrow G_{k_{3}%
}\twoheadleftarrow\cdots$ of $\mathcal{G}$, and a commutative diagram:
\[%
\begin{array}
[c]{ccccccccc}%
\pi_{1}(V_{0}) & \overset{\lambda_{1}}{\longleftarrow} & \pi_{1}(V_{1}) &
\overset{\lambda_{2}}{\longleftarrow} & \pi_{1}(V_{2}) & \overset{\lambda_{3}%
}{\longleftarrow} & \pi_{1}(V_{3}) & \overset{\lambda_{4}}{\longleftarrow} &
\cdots\\
& \overset{g_{1}}{\nwarrow}\quad\overset{f_{1}}{\swarrow} &  & \overset{g_{2}%
}{\nwarrow}\quad\overset{f_{2}}{\swarrow} &  & \overset{g_{3}}{\nwarrow}%
\quad\overset{f_{3}}{\swarrow} &  &  & \\
& G_{k_{1}} & \twoheadleftarrow &  G_{k_{2}} & \twoheadleftarrow &  G_{k_{3}}%
& \twoheadleftarrow & \cdots &
\end{array}
\]
Each $f_{i}$ is necessarily surjective, so by Lemmas \ref{gp3} and \ref{gp4}
each $\ker(f_{i})$ is the normal closure of a finite set of elements
$F_{i}\subset$ $\pi_{1}(V_{i})$. For each $i\geq1$, choose a finite collection
$\left\{  \alpha_{j}^{i}\right\}  _{j=1}^{n_{i}}$of pairwise disjoint embedded
loops in $\partial V_{i}$ representing the elements of $F_{i}$. By the
commutativity of the diagram, each $\alpha_{j}^{i}$ contracts in $V_{i-1}$.
For each $\alpha_{j}^{i}$ choose an embedded disk $D_{j}^{i}\subset
\overset{\circ}{V}_{i-1}$ with $\partial D_{j}^{i}=$ $\alpha_{j}^{i}$. Arrange
that the $D_{j}^{i}$'s are pairwise disjoint, and all intersections between
$\overset{\circ}{D_{j}^{i}}$ and $\bigcup\partial V_{k}$ are transverse.

In order to kill the kernels of the $f_{i}$'s, we would like to add to each
$V_{i}$ regular neighborhoods of the $D_{j}^{i}$'s. This would work if each
$D_{j}^{i}$ was contained in $V_{i-1}-\overset{\circ}{V}_{i}$; for then we
would be attaching a finite collection of $2$--handles to each $V_{i}$ and each
would kill the normal closure of its attaching $1$--sphere $\alpha_{j}^{i}$ in
$\pi_{1}(V_{i})$, and no more. Since this ideal situation may not be present,
we must first perform some alterations on the $V_{i}$'s.\medskip

\noindent\textbf{Claim}\qua {\sl There exists a nested cofinal sequence
$\left\{  V_{i}^{\prime}\right\}  $ of $0$--neighborhoods of
infinity which satisfy the following properties for all $i\geq1$:

\begin{itemize}
\item [\rm{(i)}]$V_{i}^{\prime}\subset V_{i}$,

\item[\rm{(ii)}] $\pi_{1}(V_{i}^{\prime})\rightarrow\pi_{1}(V_{i})$
is an isomorphism,

\item[\rm{(iii)}] $\bigcup_{j=1}^{n_{i}}\,\alpha_{j}^{i}\subset\partial
V_{i}^{\prime}$, and

\item[\rm{(iv)}] each $\alpha_{j}^{i}$bounds a 
$2$--disk in $V_{i-1}^{\prime}-\overset{\circ}{V_{i}^{\prime}}%
$.\medskip
\end{itemize}}

Roughly speaking, a $V_{q}^{\prime}$ will be constructed by removing regular
neighborhoods of the $D_{j}^{q}$'s from $V_{q}$; but in order arrange
condition (iii) and to maintain ``nestedness'', some extra care must be taken.

We already have that $\partial D_{j}^{i}=\alpha_{j}^{i}\subset\partial V_{i}$
and $\overset{\circ}{D_{j}^{i}}$ intersects finitely many $\partial V_{l}$
($l\geq i$) transversely. In addition, we would like the outermost component
of $\ D_{j}^{i}-$ $\partial V_{i}$ to lie in $V_{i-1}-\overset{\circ}{V_{i}}$.
If this is not already the case, it can easily be arranged by pushing a small
annular neighborhood of $\ \partial D_{j}^{i}$ into $V_{i-1}-\overset{\circ
}{V_{i}}$ while leaving $\partial D_{j}^{i}=\,\alpha_{j}^{i}$ fixed. Now
choose a pairwise disjoint collection $\left\{  L_{j}^{i}\right\}  $ of
regular neighborhoods of the collection $\left\{  D_{j}^{i}\right\}  $; then
for each $D_{j}^{i}$, choose a smaller regular neighborhood $N_{j}^{i}%
\subset\overset{\circ}{L_{j}^{i}}$. Between each $N_{j}^{i}$ and $L_{j}^{i}$
there exists a sequence $2N_{j}^{i},3N_{j}^{i},4N_{j}^{i},\cdots$ of regular
neighborhoods of $D_{j}^{i}$ such that
\[
N_{j}^{i}\subset\overset{\circ}{2N_{j}^{i}}\subset2N_{j}^{i}\subset
\overset{\circ}{3N_{j}^{i}}\subset3N_{j}^{i}\subset\cdots\subset L_{j}^{i}.
\]
For each $q\geq1$, let
\[
V_{q}^{\prime}=\overline{V_{q}-\bigcup\nolimits_{i=1}^{q}\left(
\bigcup\nolimits_{j=1}^{n_{i}}qN_{j}^{i}\right)  }\text{.}%
\]
Conditions (i), (iii) and (iv) are obvious, and since each $V_{i}^{\prime}$
was obtained from $V_{i\text{ }}$ by removing regular neighborhoods of
$2$--complexes, condition (ii) follows from general position.

Now, along each $\alpha_{j}^{i}$ it is possible to attach an ambiently
embedded $2$--handle $h_{j}^{i}\subset V_{i-1}^{\prime}-\overset{\circ}%
{V_{i}^{\prime}}$ to $V_{i}^{\prime}$. For each $q\geq1$, let
\[
V_{q}^{\prime\prime}=V_{q}^{\prime}\cup\left(  \bigcup\nolimits_{j=1}^{n_{q}%
}h_{j}^{q}\right)  .
\]
The naturally induced homomorphisms $f_{q}^{^{\prime\prime}}\co \pi_{1}%
(V_{q}^{\prime\prime})\rightarrow G_{k_{q}}$ are now isomorphisms.

Lastly, we must apply Lemma \ref{lemma1ngh} to each $V_{q}^{\prime\prime}$ to
create a sequence $\left\{  U_{q}\right\}  $ of generalized $1$--neighborhoods
of infinity with the same fundamental groups. To regain nestedness, we may
then have to pass to a subsequence of $\left\{  U_{q}\right\}  $ (and to the
corresponding subsequence of $\left\{  G_{k_{q}}\right\}  $) to complete the proof.
\end{proof}

The main consequences of this section are summarized by the following:

\begin{theorem}
[Generalized $1$--Neighborhoods Theorem]\label{1nghtheorem}Let $M^{n}$
($n\geq5$) be a one ended, open $n$--manifold which is inward tame at infinity. Then:

\begin{enumerate}
\item $M^{n}$ contains a neat sequence $\left\{  U_{i}\right\}  $ of
generalized $1$--neighborhoods of infinity,

\item  if $\pi_{1}(\varepsilon(M^{n}))$ is stable, we may arrange that the
$U_{i}$'s are strong $1$--neighborhoods of infinity,

\item  if $\pi_{1}(\varepsilon(M^{n}))$ is semistable, we may arrange that
each $\pi_{1}\left(  U_{i}\right)  \leftarrow\pi_{1}\left(  U_{i+1}\right)  $
is surjective, and

\item  if $\pi_{1}(\varepsilon(M^{n}))$ is perfectly semistable, we may
arrange that each $\pi_{1}\left(  U_{i}\right)  \leftarrow\pi_{1}\left(
U_{i+1}\right)  $ is surjective and has perfect kernel.
\end{enumerate}
\end{theorem}

\begin{proof}
Claim 1 is just Corollary \ref{cor1ngh}. To obtain Claim 2, observe that if\break
$\pi_{1}(\varepsilon(M^{n}))$ is pro-equivalent to $\left\{  G,id\right\}  $,
then Lemma \ref{gp3} implies that $G$ is finitely presentable. Hence we may
apply Lemma \ref{lemma1seq} to obtain the desired sequence. Claims 3 and 4
follow similarly from Lemma \ref{lemma1seq}, with the necessary algebra being
found in Lemma \ref{gp5} and Corollary \ref{gp2}.
\end{proof}

\section{Obtaining generalized $(n-3)$--neighborhoods of infinity\label{n-3}}

We now show how to obtain appropriate neat sequences of generalized $\left(
n-3\right)  $--neighborhoods of infinity. To do this, we begin with a neat
sequence $\left\{  U_{i}\right\}  $ of generalized $1$--neighborhoods of
infinity and make geometric alterations to kill $\pi_{j}(U_{i},\partial
U_{i})$ for $2\leq j\leq n-3$. These alterations will not change the
fundamental groups of the original $U_{i}$'s, hence any work accomplished by
Theorem \ref{1nghtheorem} will be preserved.

If $U_{i}$ is a generalized $1$--neighborhood of infinity and $\rho
\co \widetilde{U}_{i}\rightarrow U_{i}$ is the universal covering projection,
then $\partial\widetilde{U}_{i}=\rho^{-1}\left(  \partial U_{i}\right)  $ is
the universal cover of $\partial U_{i}$, thus, $\pi_{j}(U_{i},\partial
U_{i})\cong\pi_{j}(\widetilde{U}_{i},\partial\widetilde{U}_{i})$ for all $j$.
Moreover, if $U_{i}$ is a generalized $\left(  k-1\right)  $--neighborhood of
infinity, the Hurewicz Theorem (Theorem 7.5.4 of \cite{Sp}) implies that
$\pi_{k}(\widetilde{U}_{i},\partial\widetilde{U}_{i})$ $\cong H_{k}%
(\widetilde{U}_{i},\partial\widetilde{U}_{i})$. The last of these---the
homology in the universal cover---is usually the easiest to work with.
Throughout the remainder of this paper, the symbol ``$\sim$'' over a space
denotes a universal cover.

When calculating homology groups we prefer cellular homology. If $\left(
X^{n},Y^{n-1}\right)  $ is a manifold pair with $Y^{n-1}\subset\partial X^{n}
$, then a handle decomposition of $X^{n}$ built on $Y^{n-1}$ gives rise to a
relative CW--complex $\left(  K,Y^{n-1}\right)  \simeq\left(  X^{n}%
,Y^{n-1}\right)  $ obtained by collapsing handles onto their cores such that
each $j$--cell of $K-Y^{n-1}$ corresponds to a unique $j$--handle of $X^{n}$.
Then the cellular chain complex
\begin{equation}
0\rightarrow C_{n}\rightarrow C_{n-1}\rightarrow\cdots\rightarrow
C_{0}\rightarrow0 \tag{$\dagger$}%
\end{equation}
for $\left(  K,Y^{n-1}\right)  $, where each $C_{j}$ is generated by the
$j$--cells of $K-Y^{n-1}$, may be used to calculate the homology of $\left(
X^{n},Y^{n-1}\right)  $. We will frequently abuse terminology slightly by
referring to ($\dagger$) as the chain complex for $\left(  X^{n}%
,Y^{n-1}\right)  $ and referring to the $j$--handles of $X^{n}$ as the
generators of $C_{j}$.

If $\pi_{1}\left(  Y^{n-1}\right)  \overset{\cong}{\rightarrow}$ $\pi
_{1}\left(  X^{n}\right)  $ and we wish to calculate $H_{\ast}\left(
\widetilde{X}^{n},\widetilde{Y}^{n}\right)  $, we may use the cellular chain
complex
\begin{equation}
0\rightarrow\widetilde{C}_{n}\rightarrow\widetilde{C}_{n-1}\rightarrow
\cdots\rightarrow\widetilde{C}_{0}\rightarrow0 \tag{$\ddagger$}%
\end{equation}
of the pair $\left(  \widetilde{K},\widetilde{Y}^{n-1}\right)  $. This may be
given the structure of a $\mathbb{Z[\pi}_{1}K]$--complex, where $\widetilde
{C}_{j}$ is a free $\mathbb{Z[\pi}_{1}Y^{n-1}]$--module with one generator for
each $j$--cell of $K-Y^{n-1}$ (see Chapter I of \cite{Co} for details).
Alternatively, we will refer to ($\ddagger$) as a chain complex for $\left(
\widetilde{X}^{n},\widetilde{Y}^{n}\right)  $ where $\widetilde{C}_{j}$ has
one $\mathbb{Z[\pi}_{1}X^{n}]$--generator for each $j$--handle of $X^{n}$. The
additional algebraic structure means that each $H_{j}\left(  \widetilde{X}%
^{n},\widetilde{Y}^{n}\right)  $ may be viewed as a $\mathbb{Z[\pi}_{1}X^{n}]
$--module.

Another useful way to view ($\ddagger$) is as the chain complex for the
homology of \ $\left(  X^{n},Y^{n-1}\right)  $ with local $\mathbb{Z[\pi}%
_{1}X^{n}]$--coefficients. Then $\widetilde{C}_{j}=C_{j}\otimes\mathbb{Z[\pi
}_{1}X^{n}]$ is generated by the $j$--handles of $X^{n}$ (with preferred base
paths) and for $j>2$ the boundary map is determined by $\mathbb{Z[\pi}%
_{1}X^{n}]$--intersection numbers. In particular, if $h^{j}$ is a $j$--handle of
$X^{n}$ with attaching $\left(  j-1\right)  $--sphere $\alpha^{j-1}$, then
\[
\partial h^{j}=\sum_{s}\varepsilon\left(  \alpha^{j-1},\beta_{s}^{n-j}\right)
h_{s}^{j-1}%
\]
where $\varepsilon\left(  \alpha^{j-1},\beta_{s}^{n-j}\right)  $ denotes the
$\mathbb{Z[\pi}_{1}X^{n}]$--intersection numbers between $\alpha^{j-1}$ and the
belt sphere $\beta_{s}^{n-j}$ of a $\left(  j-1\right)  $--handle $h_{s}^{j-1}$
measured in $\partial X_{j-1}$ where $X_{j-1}=(Y^{n-1}\times\left[
0,1\right]  )\cup\left(  \text{handles of index }\leq j-1\right)  $ \ See
Chapter 6 and Appendix A of \cite{RS} for further discussion. When the chain
complex is viewed in this manner, we will still denote the corresponding
homology groups by $H_{\ast}\left(  \widetilde{X}^{n},\widetilde{Y}%
^{n}\right)  $.

The following algebraic lemma will be used each time we attempt to improve a
generalized $j$--neighborhood of infinity to a generalized $\left(  j+1\right)
$--neighborhood.

\begin{lemma}
\label{fg}Suppose $U$ is a generalized $j$--neighborhood of infinity ($j\geq1
$) in a one ended inward tame open $n$--manifold. Then $H_{j+1}\left(
\widetilde{U},\partial\widetilde{U}\right)  $ is finitely generated as a
$\mathbb{Z[\pi}_{1}U]$--module.
\end{lemma}

\begin{proof}
Fix a triangulation of $\partial U$ and for each $k\geq2$, let $K^{k}$ denote
the corresponding $k$--skeleton. Let $K^{1}$ denote the corresponding
$2$--skeleton. Note that the inclusion $K^{k}\hookrightarrow U$ induces a
$\pi_{1}$--isomorphism for all $k\geq1$. Hence, we have universal covers
$\widetilde{U}\supset\partial\widetilde{U}\supset\widetilde{K}^{k}$.

Since $H_{j}\left(  \partial\widetilde{U},\widetilde{K}^{j}\right)  =0$, the
long exact sequence for the triple $\left(  \widetilde{U},\partial
\widetilde{U},\widetilde{K}^{j}\right)  $ provides an epimorphism of
$\mathbb{Z[\pi}_{1}]$--modules $H_{j+1}\left(  \widetilde{U},\widetilde{K}%
^{j}\right)  \twoheadrightarrow H_{j+1}\left(  \widetilde{U},\partial
\widetilde{U}\right)  $. Hence, the desired conclusion will follow if we can
show that $H_{j+1}\left(  \widetilde{U},\widetilde{K}^{j}\right)  $ is
finitely generated. This will follow immediately from Theorem A of \cite{Wa1}
if we can show that $H_{i}\left(  \widetilde{U},\widetilde{K}^{j}\right)  =0$
for all $i\leq j$. Again we employ the exact sequence for $\left(
\widetilde{U},\partial\widetilde{U},\widetilde{K}^{j}\right)  $:
\[
\cdots\rightarrow H_{i}\left(  \partial\widetilde{U},\widetilde{K}^{j}\right)
\rightarrow H_{i}\left(  \widetilde{U},\widetilde{K}^{j}\right)  \rightarrow
H_{i}\left(  \widetilde{U},\partial\widetilde{U}\right)  \rightarrow\cdots
\]
Clearly the first term listed vanishes for all $i\leq j$ and, by hypothesis,
so does the third term; thus, forcing the middle term to vanish.
\end{proof}

The next lemma is the key to this section.

\begin{lemma}
\label{mainn-3}Let $M^{n}$ ($n\geq5$) be open, one ended and inward tame at
infinity and let $k\leq n-3$. Then each generalized $1$--neighborhood of
infinity $V_{0}$ contains a generalized $k$--neighborhood of infinity $U_{0}$
such that $\pi_{1}(V_{0})\overset{\cong}{\longleftarrow}\pi_{1}(U_{0})$.
\end{lemma}

\begin{proof}
If $k=1$, there is nothing to prove; otherwise we assume inductively that
$k\geq2$ and each generalized $1$--neighborhood of infinity $V$ contains a
generalized $\left(  k-1\right)  $--neighborhood of infinity $U$ such that
$\pi_{1}(V)\overset{\cong}{\longleftarrow}\pi_{1}(U)$.

Now let $V_{0}$ be a generalized $1$--neighborhood of infinity. By the
inductive hypothesis, we may assume that $V_{0}$ is already a generalized
$\left(  k-1\right)  $--neighborhood of infinity. We will show how to improve
$V_{0}$ to a $k$--neighborhood of infinity.

As we noted earlier, $\pi_{i}\left(  \widetilde{V}_{0},\partial\widetilde
{V}_{0}\right)  \cong\pi_{i}\left(  V_{0},\partial V_{0}\right)  $ for all
$i$, $H_{i}\left(  \widetilde{V}_{0},\partial\widetilde{V}_{0}\right)  =0$ for
$i\leq k-1$, and $H_{k}\left(  \widetilde{V}_{0},\partial\widetilde{V}%
_{0}\right)  \cong\pi_{k}\left(  \widetilde{V}_{0},\partial\widetilde{V}%
_{0}\right)  $. Furthermore, by Lemma \ref{fg}, $H_{k}\left(  \widetilde
{V}_{0},\partial\widetilde{V}_{0}\right)  $ is finitely generated as a
$\mathbb{Z[}\pi_{1}\left(  V_{0}\right)  ]$--module.

We break the remainder of the proof into overlapping but distinct cases:

\begin{case}
$2\leq k<\frac{n}{2}$
\end{case}

Choose a finite collection of disjoint embeddings $\left(  D_{j},\partial
D_{j}\right)  \hookrightarrow(V_{0},\partial V_{0})$ of $k$--cells representing
a generating set for $\pi_{k}\left(  V_{0},\partial V_{0}\right)  $ viewed as
a $\mathbb{Z[\pi}_{1}V_{0}]$--module. Let $Q$ be a regular neighborhood of
$\partial V_{0}\cup\left(  \bigcup D_{j}\right)  $ in $V_{0}.$ Notice that
$\pi_{1}(\partial V_{0})\rightarrow\pi_{1}(Q)$ and $\pi_{1}(Q)\rightarrow
\pi_{1}(V_{0})$ are both isomorphisms. (If $k>2$ this is obvious. If $k=2$
notice that each $\partial D_{j}$ already contracts in $\partial V_{0}$ since
$V_{0}$ is a generalized $1$--neighborhood.) Thus, $\widetilde{Q}=\rho^{-1}(Q)$
is the universal cover of $Q.$

Let $U_{0}=\overline{V_{0}-Q}$. Since the $D_{j}$'s have codimension greater
than 2, then $\pi_{1}(U_{0})\rightarrow\pi_{1}(V_{0})$ is an isomorphism. It
remains to show that $U_{0}$ is a $k$--neighborhood of infinity.

To see that $\pi_{1}(\partial U_{0})\rightarrow\pi_{1}(U_{0})$ is an
isomorphism, recall from above that $\pi_{1}(Q)\overset{\cong}{\rightarrow}%
\pi_{1}(V_{0})$. Then observe that the pair $\left(  V_{0},Q\right)  $ may be
obtained from the pair $\left(  U_{0},\partial U_{0}\right)  $ by attaching
$(n-k)$--handles (the duals of the removed handles), which has no effect on
fundamental groups.

To see that $\pi_{i}(U_{0},\partial U_{0})=0$ for $i\leq k$, we will show that
the corresponding $H_{i}\left(  \widetilde{U}_{0},\partial\widetilde{U}%
_{0}\right)  $ are trivial. By excision, it suffices to show that
$H_{i}\left(  \widetilde{V}_{0},\widetilde{Q}\right)  =0$ for $i\leq k$.

For $i<k$, the triviality of $H_{i}\left(  \widetilde{V}_{0},\widetilde
{Q}\right)  $ can be deduced from the following portion of the long exact
sequence for the triple $\left(  \widetilde{V}_{0},\widetilde{Q}%
,\partial\widetilde{V}_{0}\right)  $:
\[
\cdots\rightarrow H_{i}\left(  \widetilde{V}_{0},\partial\widetilde{V}%
_{0}\right)  \rightarrow H_{i}\left(  \widetilde{V}_{0},\widetilde{Q}\right)
\rightarrow H_{i-1}\left(  \widetilde{Q},\partial\widetilde{V}_{0}\right)
\rightarrow\cdots
\]
The first term is trivial because $V_{0}$ is a generalized $(k-1)$%
--neighborhood, and the last term is trivial (when $i-1<k)$ because
$\widetilde{Q}$ is homotopy equivalent to a space obtained by attaching
$k$--cells to $\partial\widetilde{V}_{0}$. Thus the middle term vanishes.

In dimension $k$, we use a portion of the same long exact sequence:
\[
\cdots\rightarrow H_{k}\left(  \widetilde{Q},\partial\widetilde{V}_{0}\right)
\overset{\psi}{\rightarrow}H_{k}\left(  \widetilde{V}_{0},\partial
\widetilde{V}_{0}\right)  \overset{\phi}{\rightarrow}H_{k}\left(
\widetilde{V}_{0},\widetilde{Q}\right)  \rightarrow H_{k-1}\left(
\widetilde{Q},\partial\widetilde{V}_{0}\right)  \rightarrow\cdots
\]
The last term above is again trivial for the reason cited above. Furthermore,
the map $\psi$ is surjective by the construction of $Q$; hence $\phi$ is
trivial, so $H_{k}\left(  \widetilde{V}_{0},\widetilde{Q}\right)  $ vanishes.

\begin{case}
$2<k\leq n-3$
\end{case}

The strategy in this case is similar to the above except that when $k\geq
\frac{n}{2}$ we cannot rely on general position to obtain embedded $k$--disks.
Instead we will use the tools of handle theory.

By the inductive hypothesis and the fact that $\pi_{k}\left(  V_{0},\partial
V_{0}\right)  $ is finitely generated as a $\mathbb{Z[}\pi_{1}(V_{0}%
)]$--module, we may choose a generalized $(k-1)$--neighborhood $V_{1}\subset
V_{0} $ so that, for $R=V_{0}-\overset{\circ}{V}_{1}$, the map $\pi_{k}\left(
R,\partial V_{0}\right)  \rightarrow\pi_{k}\left(  V_{0},\partial
V_{0}\right)  $ is surjective. Applying VanKampen's theorem to $V_{0}%
=R\cup_{\partial V_{1}}V_{1}$ shows that $\pi_{1}(R)\rightarrow$ $\pi
_{1}(V_{0}) $ is an isomorphism, and it follows that $\pi_{1}(\partial
V_{0})\rightarrow$ $\pi_{1}(R)$ is also an isomorphism. Hence $\rho^{-1}(R)$
is the universal cover $\widetilde{R}$ of $R$.\medskip

\noindent\textbf{Claim}\qua $H_{i}\left(  \widetilde{R},\partial\widetilde{V}%
_{0}\right)  =0$ {\sl for} $i\leq k-2$\medskip

We deduce this claim from the long exact sequence of the triple $\left(
\widetilde{V}_{0},\widetilde{R},\partial\widetilde{V}_{0}\right)  $:
\[
\cdots\rightarrow H_{i+1}\left(  \widetilde{V}_{0},\widetilde{R}\right)
\rightarrow H_{i}\left(  \widetilde{R},\partial\widetilde{V}_{0}\right)
\rightarrow H_{i}\left(  \widetilde{V}_{0},\partial\widetilde{V}_{0}\right)
\rightarrow\cdots
\]
The third term listed above is trivial for $i\leq k-1$, therefore it suffices
to show that the first term vanishes when $i\leq k-2$. Let $\widehat{V}%
_{1}=\rho^{-1}(V_{1})\subset\widetilde{V}_{0}$. Since $\pi_{1}(V_{1}%
)\rightarrow\pi_{1}(V_{0})$ needn't be an isomorphism, $\widehat{V}_{1}$
needn't be the universal cover of $V_{1}$. In fact, $\widehat{V}_{1}$ will be
connected if and only if $\pi_{1}(V_{1})\rightarrow\pi_{1}(V_{0})$ is
surjective. In general, $\widehat{V}_{1}$ has path components $\left\{
\widehat{V}_{1}^{\xi}\right\}  _{\xi\in A}$ (one for each element of
$co\ker(\pi_{1}(V_{1})\rightarrow\pi_{1}(V_{0}))$) each of which is a covering
space for $V_{1}$. See Figure 2.%
\begin{figure}
[ht!]
\begin{center}
\includegraphics[
height=3.5405in,
width=3.122in
]%
{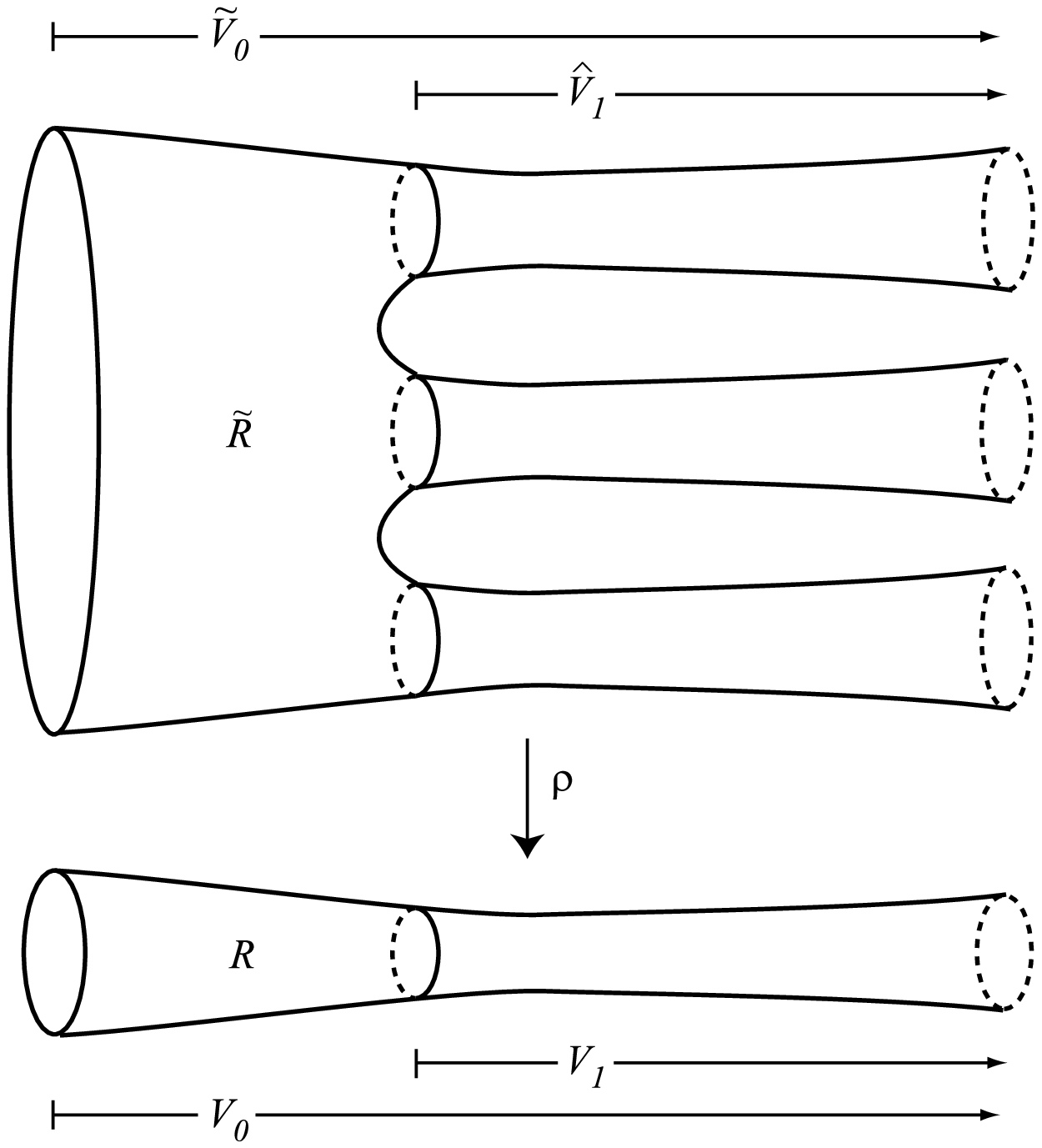}%
\caption{}%
\end{center}
\end{figure}
Moreover, $\partial\widehat{V}_{1}=\rho^{-1}(\partial V_{1})=\bigsqcup_{\xi\in
A}\partial\widehat{V}_{1}^{\xi}$, and for each $\xi$ we have $\pi_{1}%
(\partial\widehat{V}_{1}^{\xi})\overset{\cong}{\rightarrow}\pi_{1}(\widehat
{V}_{1}^{\xi})$. Thus each $\left(  \widehat{V}_{1}^{\xi},\partial\widehat
{V}_{1}^{\xi}\right)  $ is a ``covering pair'' for $\left(  V_{1},\partial
V_{1}\right)  $. It follows that $\pi_{i}\left(  \widehat{V}_{1}^{\xi
},\partial\widehat{V}_{1}^{\xi}\right)  $ is trivial for all $i\leq k-1$, so
by the Hurewicz Theorem, $H_{i}\left(  \widehat{V}_{1}^{\xi},\partial
\widehat{V}_{1}^{\xi}\right)  =0$ for all $\xi$ and for all $i\leq k-1$.
Therefore $H_{i}\left(  \widehat{V}_{1},\partial\widehat{V}_{1}\right)  =0$
for $i\leq k-1$, implying (via excision) that $H_{i+1}\left(  \widetilde
{V}_{0},\widetilde{R}\right)  $ vanishes for $i\leq k-2$, thus completing the
proof of the claim.

We now have a cobordism $\left(  R,\partial V_{0},\partial V_{1}\right)  $
with $\pi_{1}(\partial V_{0})\overset{\cong}{\rightarrow}$ $\pi_{1}(R)$ and
$H_{i}\left(  \widetilde{R},\partial\widetilde{V}_{0}\right)$ $=0$ for $i\leq
k-2$ (where $k-2\leq n-4$). By Chapter 6 of \cite{RS}, there is a handle
decomposition of $R$ built upon $\partial V_{0}$ which contains no handles of
index $\leq k-2$ and so that the existing handles have been attached in order
of increasing index. This give rise to a cellular chain complex for the pair
$\left(  \widetilde{R},\partial\widetilde{V}_{0}\right)  $ of the form:
\[
0\longrightarrow\widetilde{C}_{n}\overset{\partial_{n}}{\longrightarrow
}\widetilde{C}_{n-1}\overset{\partial_{n}}{\longrightarrow}\cdots
\overset{\partial_{k+1}}{\longrightarrow}\widetilde{C}_{k}\overset
{\partial_{k}}{\longrightarrow}\widetilde{C}_{k-1}\longrightarrow0
\]
where each $\widetilde{C}_{i}$ is a finitely generated free $\mathbb{Z[\pi
}_{1}(R)]\cong\mathbb{Z[\pi}_{1}(V_{0})]$--module with one generator for each
$i$--handle of $\left(  R,\partial V_{0}\right)  $. For $\left[  c\right]  \in$
$H_{k}\left(  \widetilde{R},\partial\widetilde{V}_{0}\right)  $ write
$c=\sum\phi_{i}e_{i}$ where each $\phi_{i}\in\mathbb{Z[\pi}_{1}(R)]$ and each
$e_{i}$ is a $k$--handle of $R$ with a preferred base path. Let $R_{k-1}\subset
R$ denote $S_{i}\cup\left(  k-1\right)  $--handles, where $S_{0}\approx\partial
V_{0}\times\left[  0,1\right]  $ is a closed collar on $\partial V_{0}$ in
$V_{0}$. We may represent $\left[  c\right]  $ with a single $k$--handle as
follows: introduce a trivial cancelling $\left(  k,k+1\right)  $--handle pair
$\left(  h^{k},h^{k+1}\right)  $ to $\partial R_{k-1}$, then do a finite
sequence of handle slides of $h^{k}$ over the other $k$--handles until $h^{k}$
is homologous to $c$. (Again see \cite{RS}.) Now, since $\partial_{k}c=0$, we
may apply the Whitney Lemma in $\partial R_{k-1}$ to move the attaching
$(k-1)$--sphere of $h^{k}$ off the belt spheres of all the $\left(  k-1\right)
$--handles.\smallskip\ 

\noindent\textbf{Note}\qua In the case $k-1=2$, the belt spheres of the
$(k-1)$--handles have codimension $2$ in $\partial R_{k-1}$, so a special case
of the Whitney Lemma (p. 72 of \cite{RS}) is needed. In particular we need to
know that the belt spheres are $\pi_{1}$--negligible in $\partial R_{k-1}$,
ie, that $\pi_{1}(\partial R_{k-1}-\left\{  \text{belt spheres}\right\}
)\overset{\cong}{\rightarrow}\allowbreak\pi_{1}(\partial R_{k-1})$. Since
$\pi_{1}(\partial V_{0})\overset{\cong}{\rightarrow}\pi_{1}(R_{k-1})$
(attaching the $2$--handles does not kill any $\pi_{1}$), this condition is
satisfied. See Lemma \ref{pi1neglemma} for the dual version of this fact.
\smallskip

We may now assume that $h^{k}$ was attached directly to $S_{0}$. By repeating
this for each element of a finite generating set for $H_{k}\left(
\widetilde{V}_{0},\partial\widetilde{V}_{0}\right)  $ we obtain a finite set
$\left\{  h_{1}^{k},\cdots,h_{t}^{k}\right\}  $ of $k$--handles attached to
$S_{0}$, so that if $Q=S_{0}\cup\left(  \bigcup h_{j}^{k}\right)  $, then
$H_{k}\left(  Q,\partial\widetilde{V}_{0}\right)  \rightarrow H_{k}\left(
\widetilde{V}_{0},\partial\widetilde{V}_{0}\right)  $ is surjective. The same
argument used in Case 1 will now show that $U_{0}=V_{0}-\overset{\circ}{Q}$ is
a generalized $k$--neighborhood of infinity.
\end{proof}

Combining Lemma \ref{mainn-3} with the Generalized $1$--Neighborhoods Theorem gives:

\begin{theorem}
[Generalized $\left(  n-3\right)  $--Neighborhoods Theorem]Let $M^{n}$
($n\geq5$) be a one ended, open $n$--manifold that is inward tame at infinity. Then

\begin{enumerate}
\item $M^{n}$ contains a neat sequence $\left\{  U_{i}\right\}  $ of
generalized $\left(  n-3\right)  $--neighborhoods of infinity,

\item  if $\pi_{1}(\varepsilon(M^{n}))$ is stable, we may arrange that the
$U_{i}$'s are strong $\left(  n-3\right)  $--neigh\-bor\-hoods of infinity,

\item  if $\pi_{1}(\varepsilon(M^{n}))$ is semistable, we may arrange that
each $\pi_{1}\left(  U_{i}\right)  \leftarrow\pi_{1}\left(  U_{i+1}\right)  $
is surjective, and

\item  if $\pi_{1}(\varepsilon(M^{n}))$ is perfectly semistable, we may
arrange that each $\pi_{1}\left(  U_{i}\right)  \leftarrow\pi_{1}\left(
U_{i+1}\right)  $ is surjective and has perfect kernel.
\end{enumerate}
\end{theorem}

\section{Obtaining generalized $(n-2)$--neighborhoods of infinity\label{n-2}}

Much like Siebenmann's original collaring theorem, the crucial step to
obtaining a pseudo-collar neighborhood of infinity is in improving generalized
$(n-3)$--neighborhoods of infinity to generalized $(n-2)$--neighborhoods of
infinity. Lemma \ref{n-2suffices} shows that, for manifolds with semistable
fundamental group systems at infinity, if we succeed our task is complete.

\begin{lemma}
\label{cobordism}Suppose $M^{n}$ $\left(  n\geq5\right)  $ contains
generalized $\left(  n-3\right)  $--neighborhoods of infinity $U_{1}\supset
U_{2}$ such that $\pi_{1}\left(  U_{1}\right)  \leftarrow\pi_{1}\left(
U_{2}\right)  $ is surjective, and let $R=U_{1}-\overset{\circ}{U}_{2}$. Then
$R$ admits a handle decomposition on $\partial U_{1}$ containing handles only
of index $(n-3)$ and $(n-2)$. Hence, $\left(  R,\partial U_{1}\right)  $ has
the homotopy type of a relative CW pair $\left(  K,\partial U_{1}\right)  $
such that $K-\partial U_{1}$ contains only $(n-3)$-- and $(n-2)$--cells.
\end{lemma}

\begin{proof}
Consider the cobordism $\left(  R,\partial U_{1},\partial U_{2}\right)  $.
Since $\pi_{1}\left(  U_{1}\right)  \twoheadleftarrow\pi_{1}\left(
U_{2}\right)  $ it is easy to check that $\pi_{1}\left(  R\right)
\twoheadleftarrow\pi_{1}\left(  \partial U_{2}\right)  $. Hence, $\pi
_{i}(R,\partial U_{2})=0$ for $i=0,1$ so we may eliminate all $0$-- and
$1$--handles from a handle decomposition of $R$ on $\partial U_{2}$. Then the
dual handle decomposition of $R$ on $\partial U_{1}$ has handles only of index
$\leq n-2$ and, by arguing as in the Claim of Lemma \ref{mainn-3}, we see that
$\pi_{i}(R,\partial U_{1})\cong H_{i}\left(  \widetilde{R},\partial
\widetilde{U}_{1}\right)  =0$ for $i\leq n-4$, so we may eliminate all handles
of index $\leq n-4$ from this handle decomposition. (In the process we
increase the numbers of $(n-3)$-- and $(n-2)$--handles.) Collapsing the
remaining handles to their cores gives us $(K,\partial U_{1})$.
\end{proof}

\begin{lemma}
\label{n-2suffices}Suppose $M^{n}$ $(n\geq5)$ contains a neat sequence
$\left\{  U_{i}\right\}  _{i=1}^{\infty}$ of $\left(  n-3\right)
$--neigh\-bor\-hoods of infinity with the property that $\pi_{1}\left(
U_{i}\right)  \leftarrow\pi_{1}\left(  U_{i+1}\right)  $ is surjective for all
$i$. Then

\begin{enumerate}
\item  each pair $\left(  U_{i},\partial U_{i}\right)  $ is homotopy
equivalent to a (probably infinite) relative CW pair $\left(  K_{i},\partial
U_{i}\right)  $ such that $K_{i}-\partial U_{i}$ \ contains only $(n-3)$-- and
$(n-2)$--cells;

\item  if some $U_{k}$ is an $\left(  n-2\right)  $--neighborhood of infinity,
then $\partial U_{k}\hookrightarrow U_{k}$ is a homotopy equivalence, ie,
$U_{k}$ is a homotopy collar.
\end{enumerate}
\end{lemma}

\begin{proof}
Roughly speaking, the first assertion is obtained by applying Lemma
\ref{cobordism} to each $R_{i}=U_{i+1}-\overset{\circ}{U}_{i}$. Since this
process is infinite, there are some technicalities to be dealt with. We refer
the reader to \cite{Si2} for details.

If $U_{k}$ is a generalized $(n-2)$--neighborhood of infinity then we already
know that $\pi_{i}(U_{k},\partial U_{k})=0$ for $i\leq n-2$. Moreover, our
first assertion guarantees that $H_{i}(\widetilde{U}_{k},\partial\widetilde
{U}_{k})$ is trivial for $i>n-2$. Thus, $\pi_{i}(\widetilde{U}_{k}%
,\partial\widetilde{U}_{k})\cong$ $\pi_{i}(U_{k},\partial U_{k})$ is trivial
for $i>n-2$, so by a theorem of Whitehead (see Section 7.6 of \cite{Sp}) $\partial
U_{k}\hookrightarrow U_{k}$ is a homotopy equivalence$.$
\end{proof}

With the end goal now clear, we begin the task of improving generalized
$(n-3)$--neighborhoods of infinity to generalized $(n-2)$--neighborhoods. Even
in the ideal situation where $\pi_{1}(\varepsilon)$ is stable and $U$ is a
strong $(n-3)$--neighborhood of infinity this may not be possible. Siebenmann
recognized that, in this ideal situation, the problem was captured by the Wall
finiteness obstruction of $U$. In our more general situation ($\pi
_{1}(\varepsilon\left(  M^{n}\right)  )$ semistable and $U$ a generalized
$(n-3)$--neighborhood of infinity) we will confront the same issue along with
some new problems caused by the lack of $\pi_{1}$--stability.

\begin{lemma}
\label{K0}Suppose $M^{n}$ $(n\geq5)$ contains a neat sequence $\left\{
U_{i}\right\}  _{i=1}^{\infty}$ of $\left(  n-3\right)  $--neigh\-bor\-hoods of
infinity with the property that $\pi_{1}\left(  U_{i}\right)
\twoheadleftarrow\pi_{1}\left(  U_{i+1}\right)  $ for all $i$. Then each
$H_{n-2}(\widetilde{U}_{i},\partial\widetilde{U}_{i})$ is a finitely generated
projective $\mathbb{Z[\pi}_{1}U_{i}]$--module. Moreover, as elements of
$\widetilde{K}_{0}\left(  \mathbb{Z[\pi}_{1}U_{i}]\right)  $, $\left[
H_{n-2}\left(  U_{i},\partial U_{i}\right)  \right]  =$ $\left(  -1\right)
^{n}\sigma(U_{i})$ where $\sigma\left(  U_{i}\right)  $ is the Wall finiteness
obstruction for $U_{i}$.
\end{lemma}

\begin{proof}
Finite generation of $H_{n-2}(\widetilde{U}_{i},\partial\widetilde{U}_{i})$
follows from Lemma \ref{fg}. For projectivity, consider the cellular chain
complex for the universal cover $\left(  \widetilde{K}_{i},\partial
\widetilde{U}_{i}\right)  $ of the CW pair $\left(  K_{i},\partial
U_{i}\right)  $ provided by assertion 1 of Lemma \ref{n-2suffices}
\[
0\rightarrow\widetilde{C}_{n-2}\overset{\partial}{\rightarrow}\widetilde
{C}_{n-3}\rightarrow0\text{.}%
\]
Triviality of $H_{n-3}(\widetilde{U}_{i},\partial\widetilde{U}_{i})$ implies
that $\partial$ is surjective, so we have a short exact sequence
\[
0\rightarrow\ker\partial\rightarrow\widetilde{C}_{n-2}\overset{\partial
}{\rightarrow}\widetilde{C}_{n-3}\rightarrow0
\]
which splits since $\widetilde{C}_{n-3}$ is a free $\mathbb{Z[\pi}_{1}U_{i}]
$--module. Thus $\widetilde{C}_{n-2}\cong\ker\partial\oplus\widetilde{C}_{n-3}%
$, so $H_{n-2}(\widetilde{U}_{i},\partial\widetilde{U}_{i})=\ker\partial$ is a
summand of a free module, and is therefore projective.

The identity\ $\left[  H_{n-2}\left(  U_{i},\partial U_{i}\right)  \right]
=\left(  -1\right)  ^{n}\sigma\left(  U_{i}\right)  $ now follows immediately
from Theorem 8 of \cite{Wa2}. An alternative argument which relies only on
\cite{Wa1} can be found in \cite{Si1}.
\end{proof}

\begin{remark}
In the above proof it is essential that $\widetilde{C}_{n-1}$ is trivial,
hence, the assumption that $\pi_{1}\left(  U_{i}\right)  \twoheadleftarrow
\pi_{1}\left(  U_{i+1}\right)  $ for all $i$ is crucial. By the Generalized
$\left(  n-3\right)  $--Neighborhoods Theorem this may be arranged whenever
$\pi_{1}(\varepsilon(M^{n}))$ is semistable (and $M^{n}$ is inward tame at infinity).
\end{remark}

\begin{lemma}
\label{1-4}Suppose $M^{n}$ ($n\geq5$) is a one ended open $n$--manifold that is
inward tame at infinity, and that $\pi_{1}(\varepsilon\left(  M^{n}\right)  )$
is semistable. Then the following are equivalent:

\begin{enumerate}
\item $M^{n}$ contains arbitrarily small clean neighborhoods of infinity
having finite homotopy types,

\item $\sigma_{\infty}\left(  M^{n}\right)  $ is trivial,

\item $M^{n}$ contains a neat sequence $\left\{  U_{i}\right\}  _{i=0}%
^{\infty}$ of generalized $\left(  n-3\right)  $--neighbor-\break
      hoods such that
$\pi_{1}\left(  U_{i}\right)  \twoheadleftarrow\pi_{1}\left(  U_{i+1}\right)
$ for all $i$ and each $H_{n-2}\left(  \widetilde{U}_{i},\partial\widetilde
{U}_{i}\right)  $ is a finitely generated stably free $\mathbb{Z}\left[
\pi_{1}U_{i}\right]  $--module,

\item $M^{n}$ contains a neat sequence $\left\{  V_{i}\right\}  _{i=0}%
^{\infty}$ of generalized $\left(  n-3\right)  $--neighbor-\break
hoods such that
$\pi_{1}\left(  V_{i}\right)  \twoheadleftarrow\pi_{1}\left(  V_{i+1}\right)
$ for all $i$ and each $H_{n-2}\left(  \widetilde{V}_{i},\partial\widetilde
{V}_{i}\right)  $ is a finitely generated free $\mathbb{Z}\left[  \pi_{1}%
V_{i}\right]  $--module.
\end{enumerate}
\end{lemma}

\begin{proof}
The equivalence of (1)--(3) follows immediately from Lemma \ref{K0} and our
earlier discussion of $\sigma_{\infty}$. Since 4)$\Longrightarrow$3) is
obvious, we need only \ show how to ``improve'' a given $U_{i}$ with stably
free $H_{n-2}\left(  \widetilde{U}_{i},\partial\widetilde{U}_{i}\right)  $ to
a generalized $\left(  n-3\right)  $--neighborhood $V_{i}$ with free
$H_{n-2}\left(  \widetilde{V}_{i},\partial\widetilde{V}_{i}\right)  $. This is
easily done by carving out finitely many trivial $\left(  n-3\right)
$--handles as described below.

Fix $i$, and let $F_{k}$ be a free $\mathbb{Z}\left[  \pi_{1}U_{i}\right]
$--module of rank $k$ so that $H_{n-2}\left(  \widetilde{U}_{i},\partial
\widetilde{U}_{i}\right)  \oplus F_{k}$ is a finitely generated free
$\mathbb{Z}\left[  \pi_{1}U_{i}\right]  $--module. Let $S_{i}\subset U_{i}$ be
a closed collar on $\partial U_{i}$ and let $\left(  h_{1}^{n-3},h_{1}%
^{n-2}\right)  ,\left(  h_{2}^{n-3},h_{2}^{n-2}\right)  ,\cdots,\left(
h_{k}^{n-3},h_{k}^{n-2}\right)  \subset$ $U_{i}-U_{i+1}$ be trivial $\left(
n-3,n-2\right)  $--handle pairs attached to $S_{i}$. Set $Q=S_{i}\cup\left(
\bigcup_{j=1}^{k}h_{j}^{n-3}\right)  $, and let $V_{i}=U_{i}-\overset{\circ
}{Q}$. VanKampen's Theorem and general position show that each of the
inclusions:  $\partial U_{i}\hookrightarrow Q$, $Q\hookrightarrow U_{i}$,
$V_{i}\hookrightarrow U_{i}$ and $\partial V_{i}\hookrightarrow V_{i}$ induce
$\pi_{1}$--isomorphisms. Thus, $V_{i}$ is a generalized $1$--neighborhood of
infinity, moreover, we have a triple $\left(  \widetilde{U}_{i},\widetilde
{Q},\partial\widetilde{U}_{i}\right)  $ of universal covers. Clearly
\[
H_{\ast}(\widetilde{Q},\partial\widetilde{U}_{i})=\left\{
\begin{array}
[c]{cc}%
0 & \text{if }\ast\neq n-3\\
F_{k} & \text{if }\ast=n-3
\end{array}
\right.  ,
\]
so for $j\leq n-3$ the long exact sequence for triples yields:
\[%
\begin{array}
[c]{ccccc}%
H_{j}\left(  \widetilde{U}_{i},\partial\widetilde{U}_{i}\right)  & \rightarrow
&  H_{j}\left(  \widetilde{U}_{i},\widetilde{Q}\right)  & \rightarrow &
H_{j-1}\left(  \widetilde{Q},\partial\widetilde{U}_{i}\right) \\
\shortparallel &  &  &  & \shortparallel\\
0 &  &  &  & 0
\end{array}
\]
Hence, $H_{j}\left(  \widetilde{U}_{i},\widetilde{Q}\right)  =0$ for $i\leq
n-3 $, and by excision, $V_{i}$ is a generalized $\left(  n-3\right)
$--neighborhood of infinity.

In dimension $n-2$ we have:
\begin{multline}%
\begin{array}
[b]{c}%
0\\
\shortparallel\\
H_{n-2}\left(  \widetilde{Q},\partial\widetilde{U}_{i}\right)
\end{array}
\rightarrow H_{n-2}\left(  \widetilde{U}_{i},\partial\widetilde{U}_{i}\right)
\rightarrow H_{n-2}\left(  \widetilde{U}_{i},\widetilde{Q}\right)
\rightarrow\\[0.2in]%
\begin{array}
[t]{c}%
H_{n-3}\left(  \widetilde{Q},\partial\widetilde{U}_{i}\right) \\
\shortparallel\\
F_{k}%
\end{array}
\rightarrow%
\begin{array}
[t]{c}%
H_{n-3}\left(  \widetilde{U}_{i},\partial\widetilde{U}_{i}\right) \\
\shortparallel\\
0
\end{array}
\nonumber
\end{multline}
Since $F_{k}$ is free this sequence splits, so
\[
H_{n-2}\left(  \widetilde{V}_{i},\partial\widetilde{V}_{i}\right)
\cong\allowbreak H_{n-2}\left(  \widetilde{U}_{i},\widetilde{Q}\right)
\cong\allowbreak H_{n-2}\left(  \widetilde{U}_{i},\partial\widetilde{U}%
_{i}\right)  \oplus F_{k}%
\]
as desired.
\end{proof}

We now begin working towards a proof of our main theorem. In order to make the
role of each hypothesis clear (and to provide additional partial results), we
begin with a minimal hypothesis and add to it only when necessary.\medskip

\noindent\textbf{Initial hypothesis}\qua {\it$M^{n}$ ($n\geq5$) is one
ended, open, and inward tame at infinity and $\pi_{1}\left(  \varepsilon
\left(  M^{n}\right)  \right)  $ is semistable.\medskip}

Then by the Generalized $\left(  n-3\right)  $--Neighborhoods Theorem we may
begin with a neat sequence $\left\{  U_{i}\right\}  _{i=0}^{\infty}$ of
generalized $\left(  n-3\right)  $--neighborhoods of infinity such that
$\pi_{1}\left(  U_{i}\right)  \twoheadleftarrow\pi_{1}\left(  U_{i+1}\right)
$ for all $i$. For each $i$, let $R_{i}=U_{i}-\overset{\circ}{U}_{i+1}$,
$\rho_{i}\co $ $\widetilde{U}_{i}\rightarrow U_{i}$ be the universal covering
projection, and $\widehat{U}_{i+1}=\rho_{i}^{-1}\left(  U_{i+1}\right)
\subset$ $\widetilde{U}_{i}$. Since each $H_{n-2}\left(  \widetilde{U}%
_{i},\partial\widetilde{U}_{i}\right)  $ is finitely generated as a
$\mathbb{Z}\left[  \pi_{1}U_{i}\right]  $--module, we may (by passing to a
subsequence and relabelling) assume that $H_{n-2}\left(  \widetilde{R}%
_{i},\partial\widetilde{U}_{i}\right)  \rightarrow H_{n-2}\left(
\widetilde{U}_{i},\partial\widetilde{U}_{i}\right)  $ is surjective for all
$i$. Consider the following portion of the long exact sequence for the triple
$\left(  \widetilde{U}_{i},\widetilde{R}_{i},\partial\widetilde{U}_{i}\right)
$.
\begin{multline*}%
\begin{array}
[b]{c}%
0\\
\shortparallel\\
H_{n-1}\left(  \widetilde{U}_{i},\widetilde{R}_{i}\right)
\end{array}
\rightarrow H_{n-2}\left(  \widetilde{R}_{i},\partial\widetilde{U}_{i}\right)
\overset{\alpha}{\rightarrow}H_{n-2}\left(  \widetilde{U}_{i},\partial
\widetilde{U}_{i}\right)  \overset{0}{\rightarrow}\\[0.2in]
H_{n-2}\left(  \widetilde{U}_{i},\widetilde{R}_{i}\right)  \overset{\gamma
}{\rightarrow}H_{n-3}\left(  \widetilde{R}_{i},\partial\widetilde{U}%
_{i}\right)  \rightarrow%
\begin{array}
[t]{c}%
H_{n-3}\left(  \widetilde{U}_{i},\partial\widetilde{U}_{i}\right)  .\\
\shortparallel\\
0
\end{array}
\end{multline*}
Triviality of the middle homomorphism follows from surjectivity of $\alpha$.
The first term in the sequence vanishes because it is isomorphic (by excision)
to $H_{n-1}\left(  \widehat{U}_{i+1},\partial\widehat{U}_{i+1}\right)  $ which
is $0$ by an application of Lemma \ref{n-2suffices}. The last term is trivial
since $U_{i}$ is a generalized $\left(  n-3\right)  $--neighborhood. After
another application of excision we obtain the following $\mathbb{Z[\pi}%
_{1}U_{i}]$--module isomorphisms:
\begin{align}
H_{n-2}\left(  \widetilde{R}_{i},\partial\widetilde{U}_{i}\right)   &  \cong
H_{n-2}\left(  \widetilde{U}_{i},\partial\widetilde{U}_{i}\right)
\tag{1}\\[0.05in]
H_{n-3}\left(  \widetilde{R}_{i},\partial\widetilde{U}_{i}\right)   &  \cong
H_{n-2}\left(  \widehat{U}_{i+1},\partial\widehat{U}_{i+1}\right)  \tag{2}%
\end{align}

By Lemma \ref{cobordism}, we may choose a handle decomposition of $R_{i}$
containing only $\left(  n-3\right)  $-- and $\left(  n-2\right)  $--handles.
Furthermore, we assume that all $\left(  n-3\right)  $--handles are attached
before any of the $\left(  n-2\right)  $--handles. Thus the homology of
$\left(  \widetilde{R}_{i},\partial\widetilde{U}_{i}\right)  $ is given by a
chain complex of the form:
\begin{equation}
0\rightarrow\widetilde{C}_{n-2}\overset{\partial}{\rightarrow}\widetilde
{C}_{n-3}\rightarrow0 \tag{3}%
\end{equation}
where $\widetilde{C}_{n-2}$ is a free $\mathbb{Z[\pi}_{1}U_{i}]$--module with
one generator for each $\left(  n-2\right)  $--handle of $R_{i}$ and
$\widetilde{C}_{n-3}$ is a free $\mathbb{Z[\pi}_{1}U_{i}]$--module with one
generator for each $\left(  n-3\right)  $--handle of $R_{i}$. From this
sequence we may extract the following short exact sequences.
\begin{align}
0  &  \rightarrow im(\partial)\rightarrow\widetilde{C}_{n-3}\rightarrow
H_{n-3}\left(  \widetilde{R}_{i},\partial\widetilde{U}_{i}\right)
\rightarrow0\tag{4}\\
0  &  \rightarrow H_{n-2}\left(  \widetilde{R}_{i},\partial\widetilde{U}%
_{i}\right)  \rightarrow\widetilde{C}_{n-2}\rightarrow im(\partial
)\rightarrow0 \tag{5}%
\end{align}
Lemma \ref{fg} (slightly modified to apply to the pair $\left(  \widetilde
{U}_{i},\widetilde{R}_{i}\right)  $) and an argument like that used in proving
Lemma \ref{K0} show that $H_{n-2}\left(  \widehat{U}_{i+1},\partial\widehat
{U}_{i+1}\right)  $ is a finitely generated projective $\mathbb{Z[\pi}%
_{1}U_{i}]$--module. Hence, by identity (2), the first of these sequences
splits. We abuse notation slightly and write
\begin{equation}
\widetilde{C}_{n-3}=im(\partial)\oplus H_{n-3}\left(  \widetilde{R}%
_{i},\partial\widetilde{U}_{i}\right)  \text{.} \tag{6}%
\end{equation}
This implies that $im(\partial)$ is also finitely generated projective, and
\begin{equation}
\left[  im(\partial)\right]  =-\left[  H_{n-3}\left(  \widetilde{R}%
_{i},\partial\widetilde{U}_{i}\right)  \right]  \in\widetilde{K}_{0}\left[
\mathbb{Z\pi}_{1}U_{i}\right]  \tag{7}%
\end{equation}
Then the second short exact sequence also splits, so we may write
\begin{equation}
\widetilde{C}_{n-2}=H_{n-2}\left(  \widetilde{R}_{i},\partial\widetilde{U}%
_{i}\right)  \oplus im(\partial)^{\prime} \tag{8}%
\end{equation}
where $im(\partial)^{\prime}$ denotes a copy of $im\left(  \partial\right)  $
lying in $\widetilde{C}_{n-2}$ (whereas $im\left(  \partial\right)  $ itself
lies in $\widetilde{C}_{n-3}$). This shows that \
\begin{equation}
\left[  im(\partial)\right]  =-\left[  H_{n-2}\left(  \widetilde{R}%
_{i},\partial\widetilde{U}_{i}\right)  \right]  \in\widetilde{K}%
_{0}(\mathbb{Z[\pi}_{1}U_{i}])\text{.} \tag{9}%
\end{equation}

\begin{remark}
Combining (1), (2), (7) and (9) shows that
\[
\left[  H_{n-2}\left(  \widetilde{U}_{i},\partial\widetilde{U}_{i}\right)
\right]  =\allowbreak\left[  H_{n-2}\left(  \widehat{U}_{i+1},\partial
\widehat{U}_{i+1}\right)  \right]  \in\allowbreak\widetilde{K}_{0}%
(\mathbb{Z[\pi}_{1}U_{i}]).
\]
In the special case that $\pi_{1}(\varepsilon(M^{n}))$ is stable and $U_{i} $
and $U_{i+1}$ are strong $\left(  n-2\right)  $--neighborhoods of infinity,
this shows that
\[
\left[  H_{n-2}\left(  \widetilde{U}_{i},\partial\widetilde{U}_{i}\right)
\right]  =\allowbreak\left[  H_{n-2}\left(  \widetilde{U}_{i+1},\partial
\widetilde{U}_{i+1}\right)  \right]  .
\]
This was one of the arguments used by Siebenmann in \cite{Si1} to show that
his end obstruction is well-defined. One can also obtain this result by using
the Sum Theorem for Wall's finiteness obstruction (see Ch. VI of \cite{Si1} or
\cite{Fe2}).
\end{remark}

Identities (6) and (8) allow us to rewrite (3) as
\begin{equation}
0\rightarrow H_{n-2}\left(  \widetilde{R}_{i},\partial\widetilde{U}%
_{i}\right)  \oplus im(\partial)^{\prime}\overset{\partial}{\rightarrow
}im(\partial)\oplus H_{n-3}\left(  \widetilde{R}_{i},\partial\widetilde{U}%
_{i}\right)  \rightarrow0 \tag{3$^{\prime}$}%
\end{equation}
where $\ker\partial=H_{n-2}\left(  \widetilde{R}_{i},\partial\widetilde{U}%
_{i}\right)  $ and $\left.  \partial\right|  _{im(\partial)^{\prime}%
}\co im(\partial)^{\prime}\overset{\cong}{\rightarrow}im(\partial)$.

We are now ready to add to our Initial Hypothesis.\medskip

\noindent\textbf{Additional Hypothesis I}\qua \emph{From now on we assume that
}$\sigma_{\infty}\left(  M^{n}\right)  =0$.\medskip

Then by Lemma \ref{1-4} we may assume that $H_{n-2}\left(  \widetilde{U}%
_{i},\partial\widetilde{U}_{i}\right)  \cong H_{n-2}\left(  \widetilde{R}%
_{i},\partial\widetilde{U}_{i}\right)  $ are finitely generated free
$\mathbb{Z[\pi}_{1}U_{i}]$--modules, and by Identity (8), that $im(\partial
)^{\prime}\cong im(\partial)$ are stably free. We may easily ``improve''
$im(\partial)^{\prime}$ and $im(\partial)$ to free $\mathbb{Z([\pi}_{1}%
U_{i}])$--modules by introducing trivial $\left(  n-3,n-2\text{ }\right)
$--handle pairs to our handle decomposition of $R_{i}$. Indeed, if we introduce
a trivial handle pair $\left(  h^{n-3},h^{n-2}\right)  $, then $im(\partial
)^{\prime}$ is increased to $im(\partial)^{\prime}\oplus\mathbb{Z[\pi}%
_{1}U_{i}]$ and $im(\partial)^{\prime}$ is increased to $im(\partial)^{\prime
}\oplus\mathbb{Z[\pi}_{1}U_{i}]$ with the new factors being generated by
$h^{h-2}$ and $h^{n-3}$, respectively. Moreover, the new boundary map
(properly restricted) is $\partial\oplus id_{\mathbb{Z[\pi}_{1}U_{i}]}$. By
doing this finitely many times we may arrange that $im(\partial)^{\prime}\cong
im(\partial)$ are free.

At this point we have a free $\mathbb{Z[\pi}_{1}U_{i}]$--module $\widetilde
{C}_{n-2}$ with a natural (geometric) basis $\left\{  h_{1}^{n-2},h_{2}%
^{n-2},\cdots,h_{r}^{n-2}\right\}  $ consisting of the $\left(  n-2\right)
$--handles of $R_{i}.$ We also have a direct sum decomposition of
$\widetilde{C}_{n-2}$ into free submodules $\widetilde{C}_{n-2}=H_{n-2}\left(
\widetilde{R}_{i},\partial\widetilde{U}_{i}\right)  \oplus im(\partial
)^{\prime}$; hence there exists another basis $\{  a_{1},\cdots
,a_{s},$\break $b_{1},\cdots,b_{r-s}\}  $ for $\widetilde{C}_{n-2}$ such that
$\left\{  a_{1},\cdots,a_{s}\right\}  $ generates $H_{n-2}\left(
\widetilde{R}_{i},\partial\widetilde{U}_{i}\right)  $ and $\left\{
b_{1},\cdots,b_{r-s}\right\}  $ generates $im(\partial)^{\prime}$. We would
like the geometry to match the algebra---in particular we would like one
subset of handles, say $\{  h_{1}^{n-2},h_{2}^{n-2},$\break $\cdots,h_{s}%
^{n-2}\}  $, to generate $H_{n-2}\left(  \widetilde{R}_{i}%
,\partial\widetilde{U}_{i}\right)  $ with the remaining handles $\{
h_{s+1}^{n-2},$\break $h_{s+2}^{n-2},\cdots,h_{r}^{n-2}\}  $ generating
$im(\partial)^{\prime}$. This may not be possible at first, but by introducing
even more trivial $\left(  n-3,n-2\text{ }\right)  $--handle pairs and then
performing handle slides, it may be accomplished. Key to the proof is the
following algebraic lemma.

\begin{lemma}
[See Lemma 5.4 of \cite{Si1}]\label{bases} Let $F$ be a finitely generated
free $\Lambda$--module with bases $\left\{  x_{1},\cdots,x_{r}\right\}  $ and
$\left\{  y_{1},\cdots,y_{r}\right\}  $ and $F^{\prime}$ be another free
module of rank $r$ with basis $\left\{  z_{1},\cdots,z_{r}\right\}  $. Then
the basis $\left\{  x_{1},\cdots,x_{r},z_{1},\cdots,z_{r}\right\}  $ of
$F\oplus F^{\prime}$ may be changed to a basis of the form $\left\{
y_{1},\cdots,y_{r},z_{1}^{\prime},\cdots,z_{r}^{\prime}\right\}  $ by a finite
sequence of elementary operations of the form $x\longmapsto x+\lambda y$.
\end{lemma}

\begin{proof}
If $A$ is the matrix of the basis $\left\{  y_{1},\cdots,y_{r}\right\}  $ in
terms of $\left\{  x_{1},\cdots,x_{r}\right\}  $, then the matrix of $\left\{
y_{1},\cdots,y_{r},z_{1},\cdots,z_{r}\right\}  $ in terms of $\left\{
x_{1},\cdots,x_{r},z_{1},\cdots,z_{r}\right\}  $ is $\left[
\begin{array}
[c]{cc}%
A & 0\\
0 & I
\end{array}
\right]  $. Now
\[
\left[
\begin{array}
[c]{cc}%
A & 0\\
0 & I
\end{array}
\right]  \cdot\left[
\begin{array}
[c]{cc}%
A^{-1} & 0\\
0 & A
\end{array}
\right]  =\left[
\begin{array}
[c]{cc}%
I & 0\\
0 & A
\end{array}
\right]  \text{,}%
\]
where
\[
\left[
\begin{array}
[c]{cc}%
A^{-1} & 0\\
0 & A
\end{array}
\right]  =\left[
\begin{array}
[c]{cc}%
I & A^{-1}\\
0 & I
\end{array}
\right]  \cdot\left[
\begin{array}
[c]{cc}%
I & 0\\
I-A & I
\end{array}
\right]  \cdot\left[
\begin{array}
[c]{cc}%
I & -I\\
0 & I
\end{array}
\right]  \cdot\left[
\begin{array}
[c]{cc}%
I & 0\\
I-A^{-1} & I
\end{array}
\right]
\]
is a product of matrices obtained by elementary moves.
\end{proof}

To apply this lemma, we introduce $r$ trivial $\left(  n-3,n-2\text{ }\right)
$--handle pairs
\[
\left\{  (k_{i}^{n-3},k_{i}^{n-2})\right\}  _{i=1}^{r}%
\]
into $R_{i}$ thus giving us a geometric basis
\[
\mathcal{B}_{1}=\left\{  h_{1}^{n-2},\cdots,h_{r}^{n-2},k_{1}^{n-2}%
,\cdots,k_{r}^{n-2}\right\}
\]
for $\widetilde{C}_{n-2}$. (In the process $\widetilde{C}_{n-2}$,
$im(\partial)^{\prime}$ and $im(\partial)$ are expanded, but\hfill\break 
$H_{n-2}\left(
\widetilde{R}_{i},\partial\widetilde{U}_{i}\right)  $ remains the same.)
According to Lemma \ref{bases} we can change $\mathcal{B}_{1}$ to a basis of
the form
\[
\mathcal{B}_{2}=\left\{  a_{1},\cdots a_{s},b_{1},\cdots b_{r-s},k_{1}%
^{\prime},\cdots,k_{r}^{\prime}\right\}
\]
using only elementary operations which may be imitated geometrically with
handle slides. Hence we arrive at a handle decomposition of $R_{i}$ for which
each element of $\mathcal{B}_{2}$ corresponds to a single $\left(  n-2\right)
$--handle, and a subset of these handles generates the submodule $H_{n-2}%
\left(  \widetilde{R}_{i},\partial\widetilde{U}_{i}\right)  $.

Given the above, we revert to our original notation in which $\widetilde
{C}_{n-2}$ is generated by $\left\{  h_{1}^{n-2},\cdots,h_{r}^{n-2}\right\}
$; assuming in addition that the subset $\left\{  h_{1}^{n-2},\cdots
,h_{s}^{n-2}\right\}  $ generates $H_{n-2}\left(  \widetilde{R}_{i}%
,\partial\widetilde{U}_{i}\right)  $.

Most of the work that remains to be done involves handle theory in the
cobordism $\left(  R_{i},\partial U_{i},\allowbreak\partial U_{i+1}\right)  $.
For convenience, we label certain subsets of $R_{i}$. Let $S_{i}\subset R_{i}$
be a closed collar on $\partial U_{i}$ , $T_{i}=S_{i}\cup\left(  \left(
n-3\right)  \text{--handles}\right)  $, and $\partial_{+}T_{i}=\partial
T_{i}-\partial U_{i}$. Then $R_{i}=T_{i}\cup\left(  h_{1}^{n-2}\cup\cdots\cup
h_{r}^{n-2}\right)  $. See Figure 3.
\begin{figure}
[ht!]
\begin{center}
\includegraphics[
trim=0.000000in 0.000000in -0.009759in 0.000000in,
height=2.9317in,
width=4.9684in
]%
{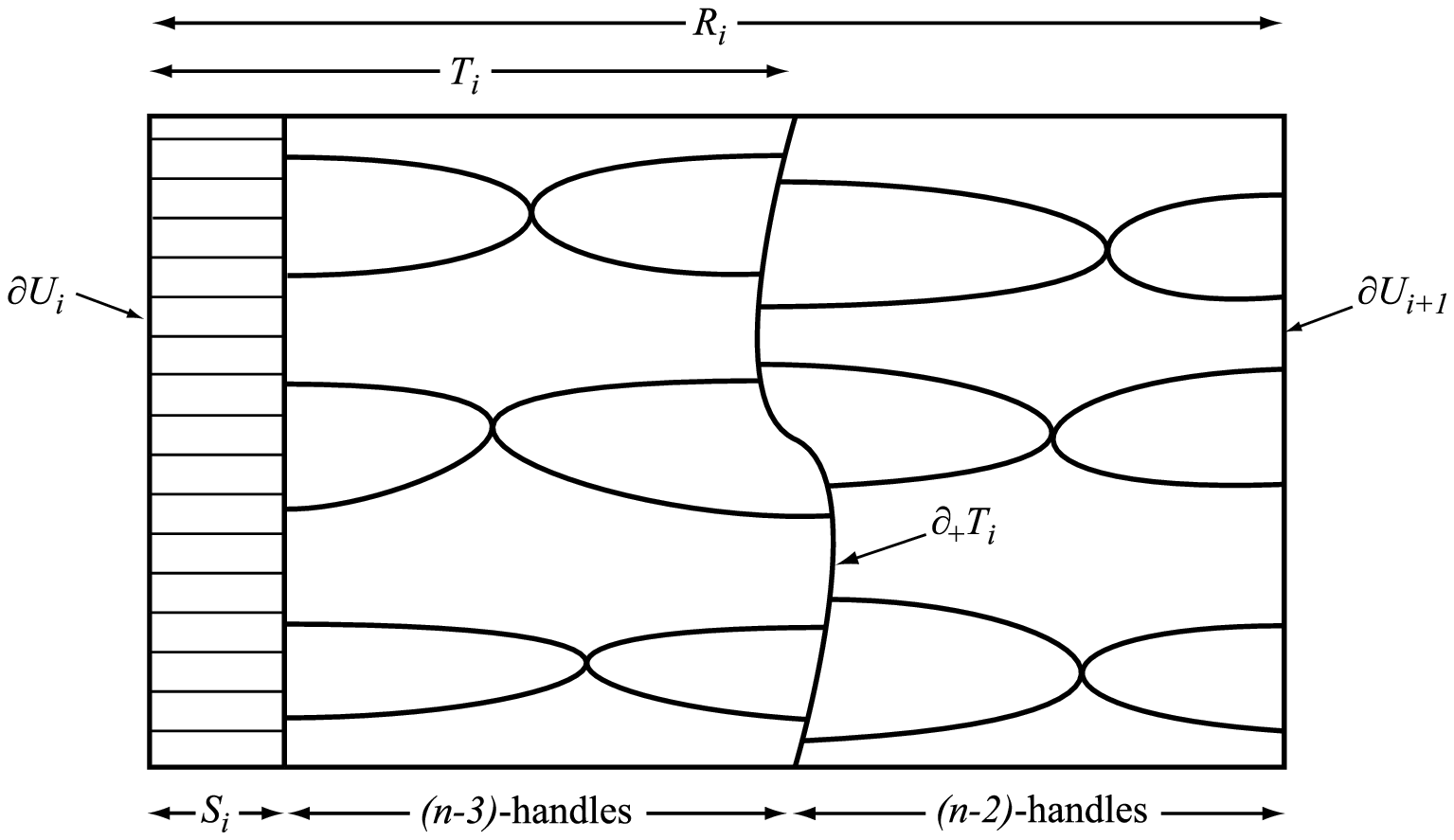}%
\caption{}%
\end{center}
\end{figure}
For each $h_{j}^{n-2}$ let $\alpha_{j}$ $\subset\partial_{+}T_{i}$ be its
attaching $\left(  n-3\right)  $--sphere, and for each $\left(  n-3\right)
$--handle $h_{k}^{n-3}$ let $\beta_{k}$ be its belt $2$--sphere.

To complete the proof, we would like to proceed as follows:

\begin{itemize}
\item  slide the handles $h_{1}^{n-2},\cdots,h_{s}^{n-2}$ (the ones which
generate $H_{n-2}\left(  \widetilde{R}_{i},\partial\widetilde{U}_{i}\right)
$) off the $\left(  n-3\right)  $--handles of $R_{i}$ so that they are attached
directly to $S_{i}$,

\item  then carve out the interiors of $h_{1}^{n-2},\cdots,h_{s}^{n-2}$ to
obtain the desired $\left(  n-2\right)  $--neighborhood.
\end{itemize}

Unfortunately, each of these steps faces a significant difficulty. For the
first step, we would like to employ the Whitney Lemma to remove $\bigcup
_{j=1}^{s}\alpha_{j}$ from $\bigcup\beta_{k}$. Since $\partial h_{j}^{n-2}=0$
for $i=1,\cdots,s$ \ the relevant $\mathbb{Z[\pi}_{1}U_{i}]$--intersection
numbers $\varepsilon\left(  \beta_{k},\alpha_{j}\right)  $ are trivial as
desired; however, since the $\alpha_{j}$'s are codimension $2 $ in
$\partial_{+}T_{i}$ we also need $\pi_{1}$--negligibility for the $\alpha_{j}%
$'s. As we will soon see, this is very unlikely.

The difficulty at the second step is similar. Assume for the moment that we
succeeded at step 1---so each $h_{j}^{n-2}$ ($j=1,\cdots s$) is attached
directly to $S_{i}$. Since the cores of the $h_{j}^{n-2}$'s have codimension
$2$ in $U_{i}$, and the $\alpha^{n-3}$'s have codimension $\ 2$ in $\partial
S_{i}$, the removal of interiors of the $h_{j}^{n-2}$'s is likely to change
the fundamental groups of these spaces---a situation we cannot tolerate at
this point in the proof.

Both of the above problems can be understood through the following easy lemma,
whose proof is left to the reader. In it, the term ``$\pi_{1}$--negligible'' is
used as follows: a subset $A$ of a space $X$ is $\pi_{1}$\emph{--negligible}
provided $\pi_{1}(X-A)\rightarrow\pi_{1}(X)$ is an isomorphism.

\begin{lemma}
\label{pi1neglemma}Suppose $\left(  W^{n},\partial_{-}W,\partial_{+}W\right)
$ is a compact cobordism ($n\geq5$) obtained by attaching $\left(  n-2\right)
$--handles $h_{1},\cdots,h_{q}$\ to a collar $C=\partial_{-}W\times\left[
0,1\right]  $. Let $\partial_{+}C$ denote $\partial_{-}W\times\left\{
1\right\}  $, and let $\alpha_{1},\cdots,\alpha_{q}\subset$ $\partial_{+}C$ be
the attaching $\left(  n-3\right)  $--spheres and $N(\alpha_{1}),\cdots
,N(\alpha_{q})\subset$ $\partial_{+}C$ the attaching tubes for the handles.
Then we have the following commutative diagram.
\[%
\begin{array}
[c]{ccc}%
\pi_{1}\left(  W^{n}\right)  & \twoheadleftarrow & \pi_{1}(\partial_{+}W)\\
\cong\uparrow &  & \cong\uparrow\\
\pi_{1}(\partial_{+}C) & \twoheadleftarrow & \pi_{1}(\partial_{+}%
C-\bigcup_{i=1}^{q}\overset{\circ}{N}(\alpha_{i}))
\end{array}
\]
Hence, the collection $\alpha_{1},\cdots,\alpha_{q}$ of attaching $\left(
n-3\right)  $--spheres is $\pi_{1}$--negligible in $\partial_{+}C$ if and only
if $\pi_{1}\left(  W^{n}\right)  \leftarrow\pi_{1}\left(  \partial
_{+}W\right)  $ is an isomorphism.
\end{lemma}

\begin{remark}
\label{siebproof}Applying this lemma to the project at hand shows that in the
special case that $\pi_{1}\left(  \partial R_{i}\right)  \leftarrow\pi
_{1}\left(  \partial U_{i+1}\right)  $ is an isomorphism for each $i$, the
program outlined above may be carried out when $n\geq6$. Hence, when $\pi
_{1}(\varepsilon\left(  M^{n}\right)  )$ is stable, we have a proof of Theorem
\ref{sieb}.
\end{remark}

Lemma \ref{pi1neglemma} shows that difficulties with fundamental groups are
unavoidable when $\pi_{1}(\varepsilon\left(  M^{n}\right)  )$ is not
stable---specifically, when $\pi_{1}\left(  U_{i}\right)  \leftarrow\pi
_{1}\left(  U_{i+1}\right)  $ has non-trivial kernel, the corresponding
attaching $\left(  n-3\right)  $--spheres will not be $\pi_{1}$--negligible.
Thus we need a new strategy for improving the $U_{i}$'s to generalized
$\left(  n-2\right)  $--neighborhoods. Instead of ``carving out'' the unwanted
$\left(  n-2\right)  $--handles in $U_{i}$, we will ``steal'' duals for these
handles from below. Our strategy is partially based on Quillen's ``plus
construction'' (see \cite{Qu} or Section 11.1 or \cite{FQ}). We will require some
additional hypotheses.\medskip

\noindent\textbf{Additional Hypothesis II}\qua $\pi_{1}\left(  \varepsilon\left(
M^{n}\right)  \right)  $\emph{\ is perfectly semistable and }$n\geq
6$\emph{.\medskip}

Then by The Generalized $\left(  n-3\right)  $--Neighborhoods Theorem we could
have chosen our original sequence $\left\{  U_{i}\right\}  _{i=0}^{\infty}$ of
generalized $\left(  n-3\right)  $--neighborhoods of infinity so that
\begin{equation}
\ker(\pi_{1}(U_{i})\twoheadleftarrow\pi_{1}(U_{i+1}))\text{\emph{is perfect
for all }}i. \tag{10}%
\end{equation}
With the exception of passing to subsequences (which is permitted by Corollary
\ref{gp2}), fundamental groups have not been changed during the current stage
of the proof, hence we may simply add Property (10) to the conditions already achieved.

Fix an $i>0$ and return to the cobordism $\left(  R_{i},\partial
U_{i},\partial U_{i+1}\right)  $ under discussion.

\begin{claim}
\label{duals}There exists a pairwise disjoint collection $\left\{  \gamma
_{j}\right\}  _{j=1}^{r}$ of embedded $2$--spheres in $\partial_{+}T_{i}$ which
are algebraic duals for the collection $\left\{  \alpha_{j}\right\}
_{j=1}^{r}$ $\ $of attaching $\left(  n-3\right)  $--spheres of the $\left(
n-2\right)  $--handles of $R_{i}$. This means that for each $0\leq j,k\leq r$,
\[
\varepsilon\left(  \alpha_{j},\gamma_{k}\right)  =\left\{
\begin{array}
[c]{c}%
1\text{\ \ if }j=k\\
\phantom{,}0\text{\ \ if }j\neq k\,.
\end{array}
\right.%
\]
\end{claim}

\noindent\textbf{Note}\qua Technically $\varepsilon\left(  \alpha_{j},\gamma
_{k}\right)  $ denotes $\mathbb{Z}\left[  \pi_{1}\left(  \partial_{+}%
T_{i}\right)  \right]  $--intersection number. Since $\pi_{1}(U_{i}%
)\overset{\cong}{\leftarrow}\pi_{1}(\partial_{+}T_{i})$ we think of it as a
$\mathbb{Z}\left[  \pi_{1}\left(  U_{i}\right)  \right]  $--intersection number.

\begin{proof}
Let $\rho\co \widetilde{U}_{i}\rightarrow U_{i}$ be the universal covering
projection. Then $\rho^{-1}(\partial_{+}T_{i})=\partial_{+}\widetilde{T}_{i}$
is the universal cover of $\partial_{+}T_{i}$. Also, $\partial_{+}%
\widetilde{T}_{i}$ $-\rho^{-1}(\bigcup_{j=1}^{r}\alpha_{j})$ covers
$\partial_{+}T_{i}-\bigcup_{j=1}^{r}\alpha_{j}$ and $\pi_{1}\left(
\partial_{+}\widetilde{T}_{i}-\rho^{-1}(\bigcup_{j=1}^{r}\alpha_{j})\right)
\cong\ker(\pi_{1}(\partial_{+}T_{i})\leftarrow\pi_{1}(\partial_{+}%
T_{i}-\bigcup_{j=1}^{r}\alpha_{j}))$ which is perfect by an application of
Lemma \ref{pi1neglemma}.

For a fixed $1\leq j_{0}\leq r$, we will show how to construct $\gamma_{j_{0}%
}$. Let $\widehat{\alpha}_{j_{0}}$ be a component of $\rho^{-1}(\alpha_{j_{0}%
})$, and let $\widehat{D}_{j_{0}}$ be a small $2$--disk in $\partial
_{+}\widetilde{T}_{i}$ intersecting $\widehat{\alpha}_{j_{0}}$ transversely in
a single point. Since $\pi_{1}\left(  \partial_{+}\widetilde{T}_{i}-\rho
^{-1}(\bigcup_{j=1}^{r}\alpha_{j})\right)  $ is perfect, then $H_{1}\left(
\partial_{+}\widetilde{T}_{i}-\rho^{-1}(\bigcup_{j=1}^{r}\alpha_{j})\right)  $
is trivial; so $\partial\widehat{D}_{j_{0}}$ bounds a surface $\widehat
{E}_{j_{0}}$ in $\partial_{+}\widetilde{T}_{i}-\rho^{-1}(\bigcup_{j=1}%
^{r}\alpha_{j})$. Let $\widehat{D}_{j_{0}}\cup\widehat{E}_{j_{0}}$ represent
an element of $H_{2}\left(  \partial_{+}\widetilde{T}_{i}\right)  $ and apply
the Hurewicz isomorphism to find a 2--sphere $\widehat{\gamma}_{j_{0}}$ in
$\partial_{+}\widetilde{T}_{i}$ representing the same element. Since they are
invariants of homology class, the $\mathbb{Z}$--intersection number of
$\widehat{\gamma}_{j_{0}}$ with $\widehat{\alpha}_{j_{0}}$ is $1$; while the
$\mathbb{Z}$--intersection number of $\widehat{\gamma}_{j_{0}}$ with any other
component of $\rho^{-1}(\bigcup_{j=1}^{r}\alpha_{j})$ is $0$. Thus, with an
appropriately chosen arc to the basepoint, the $\mathbb{Z}\left[  \pi
_{1}\left(  \partial_{+}T_{i}\right)  \right]  $--intersection numbers of
$\gamma_{j_{0}}=\rho\left(  \widehat{\gamma}_{j_{0}}\right)  $ with the
$\alpha_{j}$'s are as desired. If necessary, use general position to ensure
that $\gamma_{j_{0}}$ is embedded.
\end{proof}

By general position, we may assume that the $\gamma_{j}$'s miss the belt
$2$--spheres of each of the $\left(  n-3\right)  $--handles of $R_{i}$ ---and
hence, that they miss the $\left(  n-3\right)  $--handles altogether. Thus, the
$\gamma_{j}$'s lie in the upper boundary component $\partial_{+}S_{i}$ of the
collar $S_{i}$. The collar structure gives us a parallel copy $\gamma
_{j}^{\prime}$ $\subset\partial U_{i}$ of each $\gamma_{j}$. We would like to
arrange for each of these $\gamma_{j}^{\prime}$'s to bound a $3$--disk in
$R_{i-1}$. To make sure this is possible, we add our last additional hypotheses.\medskip

\textbf{Additional Hypotheses III}\qua $\pi_{2}(\varepsilon\left(  M^{n}\right)
)$ \emph{is semistable and} $n\geq7$.\medskip

Then, in addition to all of the above, we may assume the existence of a
diagram of the form:%
\[%
\begin{array}
[c]{ccccccccc}%
\pi_{2}(U_{0}) & \overset{\lambda_{1}}{\longleftarrow} & \pi_{2}(U_{1}) &
\overset{\lambda_{2}}{\longleftarrow} & \pi_{2}(U_{2}) & \overset{\lambda_{3}%
}{\longleftarrow} & \pi_{2}(U_{3}) & \overset{\lambda_{4}}{\longleftarrow} &
\cdots\\
& \overset{}{\nwarrow}\quad\overset{}{\swarrow} &  & \overset{}{\nwarrow}%
\quad\overset{}{\swarrow} &  & \overset{}{\nwarrow}\quad\overset{}{\swarrow} &
&  & \\
& im\lambda_{1} & \twoheadleftarrow &  im\lambda_{2} & \twoheadleftarrow &
im\lambda_{3} & \twoheadleftarrow & \cdots &
\end{array}
\]
Since each $U_{i}$ is a generalized $\left(  n-3\right)  $--neighborhood, we
have isomorphisms $\pi_{2}(\partial U_{i})\overset{\cong}{\rightarrow}\pi
_{2}(U_{i})$, $\pi_{2}(R_{i})\overset{\cong}{\rightarrow}\pi_{2}(U_{i})$ and
$\pi_{2}(\partial U_{i})\overset{\cong}{\rightarrow}\pi_{2}(R_{i})$ for all
$i\geq0$.

Assume again that $i$ has been fixed. By including the arcs to a common
basepoint (and abusing notation slightly) we view each $\gamma_{j}^{\prime} $
as representing $\left[  \gamma_{j}^{\prime}\right]  \in\pi_{2}(U_{i})$. Then,
for $1\leq j\leq r$, the above diagram guarantees the existence of a $2
$--sphere $\Omega_{j}\subset\partial U_{i+1}$ so that $\lambda_{i}\circ
\lambda_{i+1}([\Omega_{j}])=\lambda_{i}\left(  \left[  \gamma_{j}^{\prime
}\right]  ^{-1}\right)  $ in $\pi_{2}(U_{i-1})$. By general position (as
before) we may assume that $\Omega_{j}$ misses the $\left(  n-3\right)  $-- and
$\left(  n-2\right)  $--handles of $R_{i}$, and thus lies in $\left(
\partial_{+}T_{i}\right)  \cap(\partial_{+}S_{i})$ where it does not intersect
any of the $\alpha_{j}$'s. By forming the connected sum $\gamma_{j}%
\#\Omega_{j}$ (along an appropriate arc in $\partial_{+}S_{i}$) we obtain a
new $2$--sphere in $\partial_{+}T_{i}$ with $\varepsilon\left(  \alpha
_{k},\gamma_{j}\#\Omega_{j}\right)  =\varepsilon\left(  \alpha_{k},\gamma
_{j}\right)  $ for all $k$, and the additional property that its parallel copy
$(\gamma_{j}\#\Omega_{j})^{\prime}$ in $\partial U_{i}$ contracts in $U_{i-1}%
$. Note then that $(\gamma_{j}\#\Omega_{j})^{\prime}$ may be contracted in
$R_{i-1\text{.}}$

In order to simplify notation, we replace each $\gamma_{j}$ with $\gamma
_{j}\#\Omega_{j}$ and assume that, in addition to the properties of \ Claim
\ref{duals}, we have chosen the $\gamma_{j}$'s to satisfy:
\[
\emph{Each}\text{ }\gamma_{j}\text{\emph{ is contained in }}\partial_{+}%
S_{i}\text{\emph{ and its parallel copy }}\gamma_{j}^{\prime}\subset\partial
U_{i}\text{\emph{ contracts in }}R_{i-1}.
\]
By general position (here we use $n\geq7$), we may select a pairwise disjoint
collection $\left\{  D_{j}\right\}  _{j=1}^{s}$ of properly embedded $3$--disks
in $R_{i-1}$ with $\partial D_{j}=\gamma_{j}^{\prime}$ for each $1\leq j\leq
s$.\medskip

\noindent\textbf{Note}\qua We have selected bounding disks only for the duals
to\ the attaching spheres $\alpha_{1},\cdots,\alpha_{s}$ of the handles
$h_{1}^{n-2},\cdots,h_{s}^{n-2}$ which generate $H_{n-2}(U_{i},\partial
U_{i})$.\medskip

Now let $Q_{i}$ be a regular neighborhood in $R_{i-1}$ of $\partial U_{i}%
\cup\left(  \bigcup_{j=1}^{s}D_{j}\right)  $ and let $V_{i}=Q_{i}\cup U_{i}$.
See Figure 4.
\begin{figure}
[ht!]
\begin{center}
\includegraphics[
height=2.8003in,
width=4.9606in
]%
{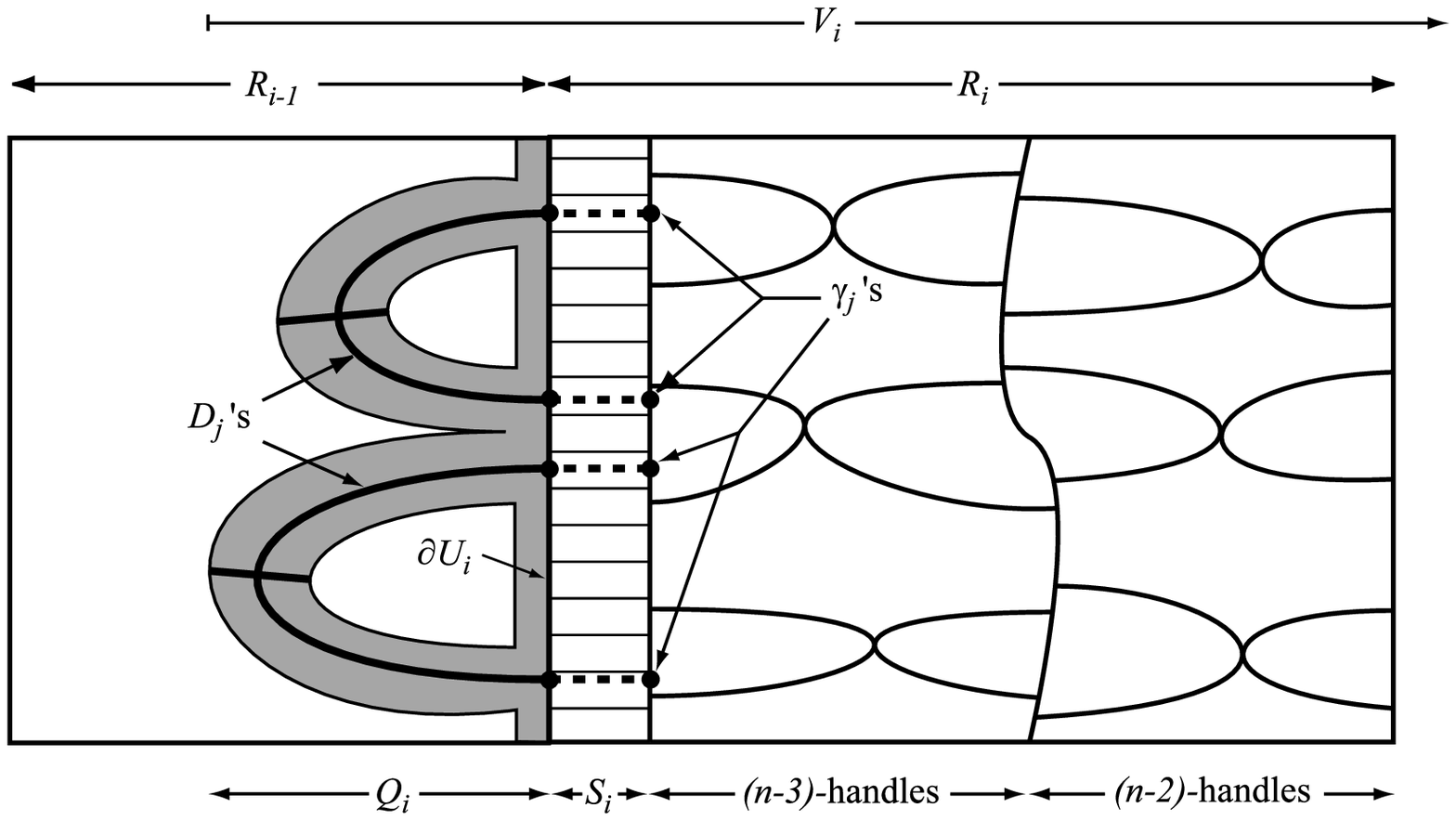}%
\caption{}%
\end{center}
\end{figure}
Our proof of the Main Existence Theorem will be complete when we prove the following.

\begin{claim}
$V_{i}$ is a homotopy collar.
\end{claim}

Notice that $\pi_{1}\left(  U_{i}\right)  \rightarrow\pi_{1}\left(
V_{i}\right)  $ is an isomorphism and $U_{i-1}\supset V_{i}\supset U_{i+1}$;
so $V_{i}$ may be substituted for $U_{i}$ as part of a $\pi_{1}$--surjective
system of neighborhoods of infinity. Hence, by Lemma \ref{n-2suffices}, it
suffices to show that $V_{i}$ is a generalized $\left(  n-2\right)
$--neighborhood of infinity.

Consider the cobordism $\left(  Q_{i},\partial V_{i},\partial U_{i}\right)  $.
Since $Q_{i}$ may be obtained by attaching $3$--handles (one for each $D_{j}$)
to a collar on $\partial U_{i}$, it may also be constructed by attaching
$\left(  n-3\right)  $--handles to a collar on $\partial V_{i}$. In this case
the $2$--spheres $\gamma_{1},\cdots,\gamma_{s}$ become the belt spheres of the
$\left(  n-3\right)  $--handles, which we label as $k_{1}^{n-3},\cdots
,k_{s}^{n-3}$, respectively. Since we already know that $U_{i}$ admits an
infinite handle decomposition containing only $\left(  n-3\right)  $-- and
$\left(  n-2\right)  $--handles, this shows that $V_{i}$ also admits a handle
decomposition with handles only of these indices. It follows from general
position that $V_{i}$ is a generalized $1$--neighborhood and from the usual
argument that $\pi_{k}(V_{i},\partial V_{i})\cong\allowbreak\pi_{k}%
(\widetilde{V}_{i},\partial\widetilde{V}_{i})\cong\allowbreak H_{k}%
(\widetilde{V}_{i},\partial\widetilde{V}_{i})=0$ for $k\leq n-4$. Hence, it
remains only to show that $H_{n-3}(\widetilde{V}_{i},\partial\widetilde{V}%
_{i})=0=H_{n-2}(\widetilde{V}_{i},\partial\widetilde{V}_{i})$. To do this,
begin with an infinite handle decomposition of $\left(  U_{i},\partial
U_{i}\right)  $ which has only $\left(  n-3\right)  $-- and $\left(
n-2\right)  $--handles, and which contains the handle decomposition of $R_{i}$
used above. Let
\begin{equation}
0\rightarrow\widetilde{D}_{n-2}\overset{\partial}{\rightarrow}\widetilde
{D}_{n-3}\rightarrow0 \tag{11}%
\end{equation}
be the associated chain complex for the $\mathbb{Z[\pi}_{1}(U_{i})]$--homology
of $\left(  U_{i},\partial U_{i}\right)  $. Then $\partial$ is surjective and
$\widetilde{D}_{n-2}=\ker\partial\oplus\widetilde{D}_{n-3}^{\prime}$, where
$\left.  \partial\right|  _{\widetilde{D}_{n-3}^{\prime}}$ is an isomorphism,
and $\ker\partial=H_{n-2}(\widetilde{U}_{i},\partial\widetilde{U}%
_{i})=\left\langle h_{1}^{n-2},\cdots,h_{s}^{n-2}\right\rangle $. Hence, (11)
may be rewritten as:
\begin{equation}
0\rightarrow\left\langle h_{1}^{n-2},\cdots,h_{s}^{n-2}\right\rangle
\oplus\widetilde{D}_{n-3}^{\prime}\overset{0\oplus\left.  \partial\right|
_{\widetilde{D}_{n-3}^{\prime}}}{\longrightarrow}\widetilde{D}_{n-3}%
\rightarrow0 \tag{11$^{\prime}$}%
\end{equation}

Our preferred handle decomposition of $(V_{i},\partial V_{i})$ is obtained by
inserting the $\left(  n-3\right)  $--handles $k_{1}^{n-3},\cdots,k_{s}^{n-3}$
beneath our handle decomposition of $\left(  U_{i},\partial U_{i}\right)  $.
Hence, the corresponding chain complex for $(V_{i},\partial V_{i})$ has the
form
\[
0\rightarrow\left\langle h_{1}^{n-2},\cdots,h_{s}^{n-2}\right\rangle
\oplus\widetilde{D}_{n-3}^{\prime}\overset{\partial_{1}\oplus\partial_{2}%
}{\longrightarrow}\left\langle k_{1}^{n-3},\cdots,k_{s}^{n-3}\right\rangle
\oplus\widetilde{D}_{n-3}\rightarrow0.
\]
In the usual way, the image of an $\left(  n-2\right)  $--handle under the
boundary map is determined by the $\mathbb{Z\pi}_{1}$--intersection numbers of
its attaching $\left(  n-3\right)  $--sphere with the belt $2$--spheres of the
various $\left(  n-3\right)  $--handles. Thus it is easy to see that $\partial
h_{j}^{n-2}=\left(  k_{j}^{n-3},0\right)  \in\left\langle k_{1}^{n-3}%
,\cdots,k_{s}^{n-3}\right\rangle \oplus\widetilde{D}_{n-3}$ for each
$h_{j}^{n-2}$ ($1\leq j\leq s$); and the map $\partial_{2}\co \widetilde{D}%
_{n-3}^{\prime}\rightarrow\left\langle k_{1}^{n-3},\cdots,k_{s}^{n-3}%
\right\rangle \oplus\widetilde{D}_{n-3}$ is of the form $\left(
\lambda,\left.  \partial\right|  _{\widetilde{D}_{n-3}^{\prime}}\right)  $
where $\lambda$ is unimportant to us and $\left.  \partial\right|
_{\widetilde{D}_{n-3}^{\prime}}$is the isomorphism from (11$^{\prime}$). It is
now easy to check that $\partial_{1}\oplus\partial_{2}$ is an isomorphism; and
thus, $H_{n-3}(\widetilde{V}_{i},\partial\widetilde{V}_{i}) $ and
$H_{n-2}(\widetilde{V}_{i},\partial\widetilde{V}_{i})$ are trivial.\medskip

\noindent\textbf{Note}\qua For those who wish to avoid the technical issues
involved with infinite handle decompositions, an alternative proof that
$V_{i}$ is a generalized $\left(  n-2\right)  $--neighborhood may be obtained
by analyzing the long exact sequence for the triple $\left(  V_{i},Q_{i}\cup
R_{i},\partial V_{i}\right)  $. The work involved is similar, but the key
calculations are now shifted to the compact pair $\left(  Q_{i}\cup
R_{i},\partial V_{i}\right)  $.

\section{Questions}

The results and examples discussed in this paper raise a number of natural
questions. We conclude this paper by highlighting a few of them.

The most obvious question is whether Conditions 1--3 of Theorem \ref{met} are
sufficient to imply pseudo-collar\-abil\-ity. Other possible improvements to
Theorem \ref{met} involve Condition 2. For example, it seems reasonable to
hope that the assumption of ``perfect semistability'' can be weakened to just
``semistability''. Note that Condition 4 and ``perfectness'' were not used
until very near the end of the proof Theorem \ref{met}.

Unlike the conditions just mentioned, the assumption that $\pi_{1}\left(
\varepsilon(M^{n})\right)  $ is semi\-stable is firmly embedded in the proof
of Theorem \ref{met}. In particular, nearly all \ of the work done in Section
\ref{n-2} depends on this assumption. However, we do not know an example of an
open manifold that is inward tame at infinity which is not $\pi_{1}%
$--semistable at infinity. Therefore we ask: Is every one ended open manifold
that is inward tame at infinity also $\pi_{1}$--semistable at infinity? Must
$\pi_{1}$ be perfectly semistable at infinity?

Lastly, we direct attention towards universal covers of closed aspherical
manifolds. As we noted in Section \ref{pseudo}, these provide some of the most
interesting examples of pseudo-collarable manifolds. Hence, we ask whether the
universal cover of a closed aspherical manifold is always pseudo-collarable.
Since very little is known in general about the ends of universal covers of
closed aspherical manifolds, one should begin by investigating whether these
examples must satisfy \emph{any} of the conditions in the statement of Theorem
\ref{met}.

\end{document}